\documentclass[12pt]{article}
\usepackage[english]{babel}
\usepackage{amsmath}
\usepackage{amssymb}
\usepackage{amsfonts}
\usepackage{amsthm}
\usepackage{graphics,graphicx}
\usepackage{color}
\usepackage{epstopdf}

\theoremstyle{plain}

\newtheorem*{proposition*}{Proposition}

\newtheorem*{corollary*}{Corollary}
\newtheorem*{Corollary*}{Important corollary}

\theoremstyle{definition}

\theoremstyle{remark}

\newtheorem*{remark*}{Remark}

\newtheorem*{example*}{Example}

\newtheorem*{examples*}{Examples}

\numberwithin{equation}{section}
\numberwithin{lemma}{section}
\numberwithin{theorem}{section}
\numberwithin{hypothesis}{section}
\numberwithin{definition}{section}
\numberwithin{example}{section}
\numberwithin{corollary}{section}
\numberwithin{remark}{section}

 \def\EE{\vbox {\hbox to 8.9pt {I\hskip-2.1pt E\hfil}}}
 \def\e{{\rm e}}
 \def\tang{\hbox{\rm tang}}
\begin{document}
\begin{center}
{\large A.Ya. KHINTCHINE's WORK
IN PROBABILITY THEORY}

\vspace{0.5mm}
{\bf Sergei Rogosin$^{1}$, Francesco Mainardi$^2$}

 \textrm{$^{1}$ Department of Economics,
 Belarusian State University
\protect\\
Nezavisimosti ave 4,  220030 Minsk, Belarus\\
e-mail: rogosin@bsu.by   $\quad$ (Corresponding author)}\\
\textrm{$^{2}$ Department of Physics and Astronomy,
Bologna University and INFN
\protect\\
Via Irnerio 46, I-40126 Bologna,  Italy\\
e-mail: francesco.mainardi@bo.infn.it}
\end{center}

\vskip -2.5truecm
%%%%%
\begin{abstract}
\noindent The paper is devoted to the contribution in the Probability Theory
of the well-known Soviet mathematician Alexander Yakovlevich Khintchine (1894--1959).
Several of his results are described, in particular those fundamental results on the infinitely divisible distributions. 
Attention is paid also to his interaction with Paul L\'evy. The content of the paper
is related to our joint book {The Legacy of A.Ya. Khintchine's Work in Probability Theory} published in 2010
by Cambridge Scientific Publishers. It  is published  in 
{ Notices of the International Congress of Chinese Mathematicians (ICCM)}, International Press, Boston, Vol. 5, No 2, pp. 60--75  (December 2017).
\\
DOI: 10.4310/ICCM.2017.v5.n2.a6
%% http://intlpress.com/site/pub/pages/journals/items/iccm/_home/_main/

%\vspace{1mm}

%\vspace{1mm}

\end{abstract}

%\noindent $^{*}${\footnotesize{Corresponding author.}}

%

\date{\today}

\vspace{0.5mm}
%{\bf CONTENTS}

{\bf 1. Introduction}$\ldots\ldots\ldots\ldots\ldots\ldots\ldots\ldots\ldots\ldots\ldots\ldots\ldots\ldots\ldots\ldots . .$ 2

{\bf 2. Short biography of Alexander Yakovlevich Khintchine}$\ldots$ 3

{\bf 3. First papers in Probability: 1924--1936}$\ldots\ldots\ldots\ldots\ldots\ldots . . . $ 5

{\bf 4. The interaction with Paul L\'evy}$\ldots\ldots\ldots\ldots\ldots\ldots\ldots\ldots\ldots . .$10

{\bf 5. Infinitely divisible distributions}$\ldots\ldots\ldots\ldots\ldots\ldots\ldots\ldots\ldots . $ 12

{\bf 6. Khintchine's book on the distribution of
the sum

\hspace{4mm} of independent random variables}$\ldots\ldots\ldots\ldots\ldots\ldots\ldots\ldots\ldots $ 14

{\bf 7. Teaching Probability Theory and Analysis}$\ldots\ldots\ldots\ldots\ldots .$ 14

\hspace{5mm}{7.1 Pedagogical credo}$\ldots\ldots\ldots\ldots\ldots\ldots\ldots\ldots\ldots\ldots\ldots\ldots\ldots\ldots .$ 15

\hspace{5mm}{7.2 Mathematics in secondary schools}$\ldots\ldots\ldots\ldots\ldots\ldots\ldots\ldots\ldots   $  $\,$17

\hspace{5mm}{7.3 Teaching mathematical analysis}$\ldots\ldots\ldots\ldots\ldots\ldots\ldots\ldots\ldots\ldots$22

\hspace{5mm}{7.4 Teaching probability theory}$\ldots\ldots\ldots\ldots\ldots\ldots\ldots\ldots\ldots\ldots\ldots . $23

{\bf Bibliography on A.Ya. Khintchine (A.I. Khinchin})
$\ldots\ldots\ldots .$26

{\bf References}$\ldots\ldots\ldots\ldots\ldots\ldots\ldots\ldots\ldots\ldots\ldots\ldots\ldots\ldots\ldots\ldots\ldots\ldots ..$41

%%%%%%%%%%

\newpage 
%%%%%%%%
\section{Introduction}

This paper deals with the work of the well-known Russian
mathematician, Alexander Yakovlevich Khintchine\footnote{Another
transliteration of his name in English is Alexander
Iacovlevich Khinchin or Hintchine.}  (1894--1959),
one of the creator of the Russian (Soviet)
mathematical school on Probability Theory. 

An extended version of the paper which includes
translation from Russian, French, German and Italian into English several articles by A.Ya.Khintchine
(as well as his main monograph) is presented in the book \cite{RogMai_Legacy}:

\noindent S.V. Rogosin, F. Mainardi, {\it The Legacy of A.Ya. Khintchine’s
Work in Probability Theory}, Cambridge Scientific Publishers (2010).

The idea to write this paper arises from two main considerations. 

First of all, Probability Theory has been developed in 1920-1930s so rapidly that
now some of the ideas and results of this period are rediscovering. It
is very interesting from different points of view to understand how this
branch of mathematics became one of the most important and applicable
mathematical discipline. In order to see it one can turn to the most
influential results concerning the main concepts of the Probability Theory. 

Second, not all results of the pioneers of this development
are known to  new generations of mathematicians. 

Among the contributors to the ``great jump" of the Probability in the above
mentioned period special role belongs to Alexandr Yakovlevich Khintchine.
He had obtained outstanding results, which form a very clear and rigorous
style of handling important deep problems. He creates together with
A. N. Kolmogorov the well-known school of Probability Theory at Moscow
University.

In connection with the above said motivation for our paper, it is especially important to describe the works by Khintchine because of the following reasons:

- Several important results by Khintchine are forgotten and later re-discovered.

- A number of results were published in inaccessible places and not in
English.

- The concrete and clear style of Khintchine's work  can help the readers to
understand  better some recent results.

%%%%%%%%%
\newpage
%%%%%%%%

The paper has the similar structure  to that of  the above mentioned book \cite{RogMai_Legacy}: S.V. Rogosin, F. Mainardi, {\it The Legacy of A.Ya. Khintchine’s
Work in Probability Theory}, Cambridge Scientific Publishers (2010). Section 5, describing the results by A.Ya. Khintchine on infinitely divisible distribution, is highly related to our paper \cite{FMSR_M2EF}:
Mainardi, F. and  Rogosin, S., The origin of infinitely divisible
distributions: from de Finetti's problem  to L\'evy-Khintchine formula,
 {\it  Mathematical Methods for Economics and Finance}, {\bf 1} (2006) 37--55.
 [E-print {\tt http://arXiv.org/math/arXiv:0801.1910}]

\vspace{2mm}
\section{Short biography of Alexander Yakovlevich Khintchine}

A.Ya. Khintchine was born on July 19, 1894 in the
village Kondrovo of the Kaluga region, about one and a half
hundred km southwest of Moscow. From 1911 to 1916 he was a student
of the Physical-Mathematical faculty of the Moscow State University (MSU).
All his scientific life was in deeply connected with this University.

In the period of study at the University and in the first years of his
research career Khintchine was under a strong influence of the
ideas and personality of N.N.~Luzin. It is known that
A.Ya.~Khintchine presented his first result at a meeting of the
student mathematical club in November 1914 (see, {\it e.g.},\,
his obituary by  Gnedenko (1961) \cite{Gnedenko_obituary}).

The mathematical talent of this young student was noticed at the
University by his teachers. After graduation at MSU A.Ya.~Khintchine was
recommended for preparation to the professorship.
His teaching career  started in 1918 at Moscow Women's Polytechnical Institute.
One year later he was invited to the Ivanovo-Voznesensk (now Ivanovo) Polytechnical Institute,
and soon after  he became the dean of Physical-Mathematical Faculty of newly founded
Ivanovo-Voznesensk Pedagogical Institute.
In 1922 the Research Institute on Mathematics and Mechanics was organized at the Moscow State
University.
A.Ya.~Khintchine was invited to this Institute as a researcher.
During a certain  period he combined his research work in
Moscow with lecturing at Ivanovo-Voznesensk.
Finally, in 1927 he  got the professorship at the Moscow State University.
%%%%%%%%%
\newpage 
%%%%%%%%
\vskip 0.5truecm
%\begin{figure}[ht]
	\begin{center}
	\includegraphics[width=0.7\textwidth]{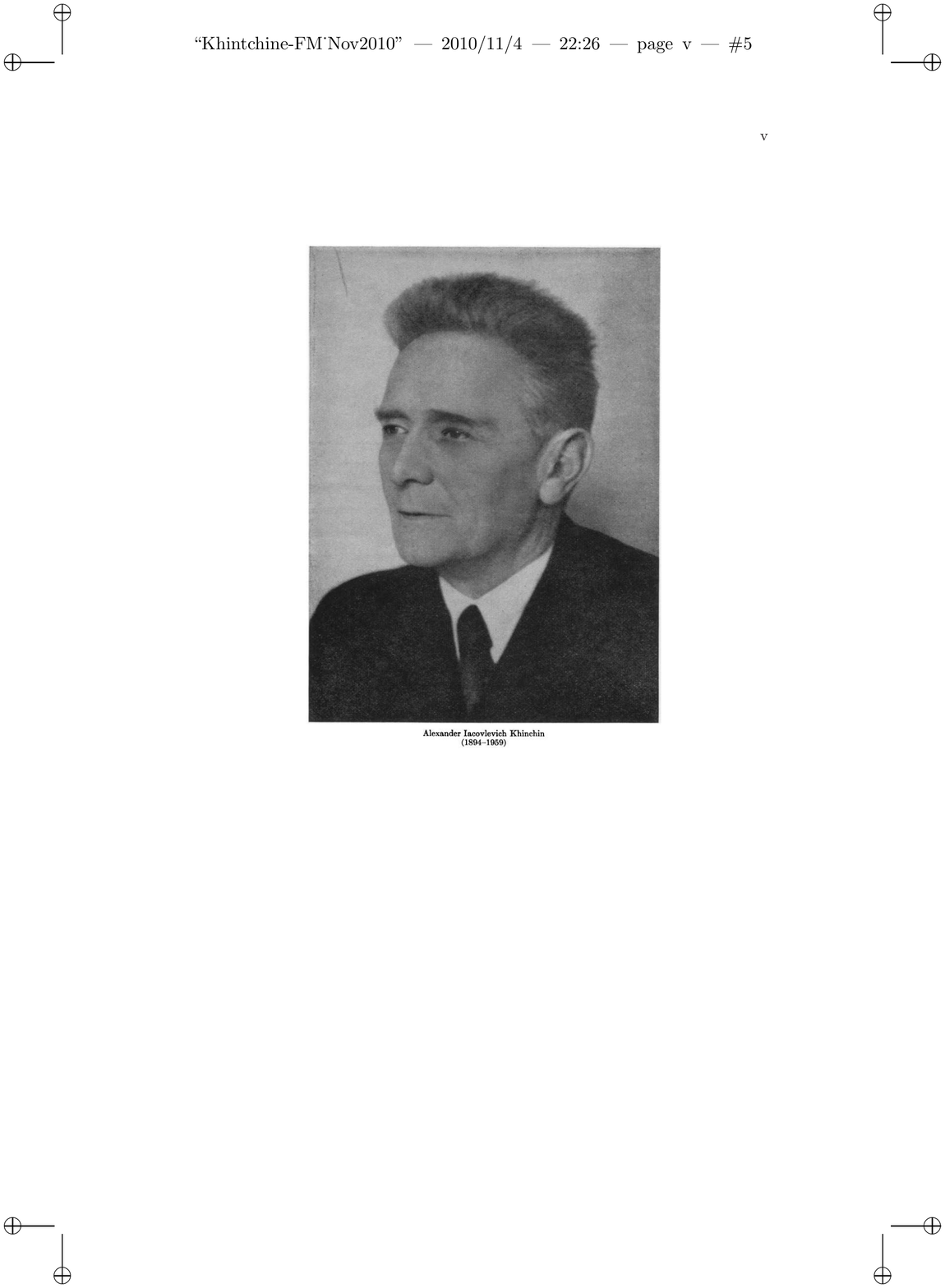}
	\end{center}

After his first significant publications A.~Ya.Khintchine became
known to the European probabilistic community. In 1928 he spent a
couple of weeks at the University in  G\"{o}ttingen,
 one of the most important mathematical centers
of the beginning of XX century.
Here he  prepared  at least two papers, published in 1929
in Mathematische Zeitschrift and in Mathematische Annalen ([G46]
and [G48], respectively).\footnote{Throughout our paper the
citation of the type [G~$\ldots$] means the corresponding paper in
the list of Khintchine publications, presented by B.~V.~Gnedenko
in his 1961 article \cite{Gnedenko_obituary}. We have included the revised
version of this list at the end of our paper.}
%%%%%%%%%%5
\newpage
%%%%%%%%%%%

Khintchine was a member of the Soviet delegation at the
International Congress of Mathematicians held in Bologna (Italy)
from 3 to 10 September 1928 (see, {\it e.g.}, Mainardi and Rogosin  (2006) \cite{FMSR_M2EF}).
%At the
%Bologna Congress de Finetti had  the  occasion  to meet L\'evy and
%Khintchine  (who were included in the French and Russian
%delegations, respectively) but we are not informed about their
%interaction.
The  Russian delegation was represented by   27  scientists
including  some prominent researchers like
%% B. Hostinsk\'y  and F.M. Urban (Brno)
%% A. Lomnicki (Poland)
S. Bernstein (Karkhov),
%% S. Timoshenko (Michigan, USA)
A.Ya.~Khintchine (Moscow),  V. Romanovsky (Tashkent) and E.
Slutsky (Moscow). We note, however, that Khintchine did not
present any communication so that he  did not publish a paper in the
Proceedings of the Congress (which appeared in 1929-1932).

Since 1927 all A.Ya.~Khintchine's further scientific and
teaching activity was connected with Moscow State University. He
was the  head of the chair of the Probability Theory, then the
head of the chair of Mathematical Analysis, director of the
Research Institute of Mathematics and Mechanics at MSU. He
passed away on November 18, 1959 after a long heavy illness.

\section{First papers in Probability: 1924--1936}

The first papers by A.~Ya.~Khintchine were appeared in 1924.
In order to understand the role of these articles
it is necessary to describe the state of Probability Theory in those years. One can recall
a critical review by R. von Mises who summed up of the situation in the following words:
``To-day, probability theory is not a mathematical science'' (see, {\it e.g.}, Cramer (1962) \cite{Cramer_obituary}).

There was no satisfactory definition of mathematical probability,
and the conceptual foundations of the subject were completely obscure.
Moreover, with few exceptions, mainly belonging to the French and Russian schools,
writers on probability did not seem aware of the standards of rigor which
in other mathematical fields, were regarded as obvious.

Already in the middle of 1920s  were appeared another
estimate of Probability Theory as a branch of mathematical
science. 

In his monograph ``Fundamental Laws of Probability''
([G~35]), published in 1927, Khintchine wrote:
``Up to the recent years in Europe the dominated opinion on Probability Theory was as on
science which is important and useful, but cannot stated and
solved serious problems. Anyway the work of Russian mathematicians
(in particular, results by P.~L.~Chebyshev, A.~M.~Lyapunov,
A.~A.~Markov) show us that it is not a correct point of view.
%%%%%%%%%%
\newpage
%%%%%%%%%
The Probability Theory has an integral method deeply connected with
the methods of modern theory of functions\footnote{Khintchine had
in mind first of all the results on measure theory and different
generalizations of the integral, which were very popular when hos book 
appeared.}, and thus the most of the
recent ideas appeared in the Mathematical Analysis have a fruitful
application in the Probability Theory.''

This optimistic opinion by A.~Ya.~Khintchine has got an evident
justification in the next few decades.

At the end of the 1930s, the picture has been radically changed.
Mathematical probability theory was firmly established on an
axiomatic foundation. It became a purely mathematical discipline, with
problems and methods of its own, conforming the current standards
of mathematical rigorism, and entering into fruitful relations with
other branches of mathematics. At the same time, the fields of
applications of mathematical probability were steadily and rapidly
growing in number and importance.
It is true that nowadays there are still %(and this may apply even today, in 1962)
some ``pure'' mathematicians who tend to look down on the ``applied'' science of probability.
But this attitude is expected to  disappear within a generation.

The tremendous development of Probability Theory, which thus took place in the twenty
years from 1920s to 1940s was, no doubt, a joint effect of the
efforts of a  number of mathematicians and statisticians.
However, it does not seem unlikely that future historians will
ascribe its development, as far as the mathematical side of the
subject is concerned, above all to the creative powers of four
scientists (in alphabetic order): B. de Finetti, A.Ya. Khintchine, A.N. Kolmogorov,
and P. L\'evy. In fact, it
may be said that the real turning point came with the publications
of the following works:

\noindent
- P. L\'evy,
  {\it Calcul des Probabilit\'es},
  Gauthier-Villars, Paris (1925), pp. viii+350.
  %%  [Reprinted in 2003  by
  %% Editions Jaques Gabay (Les Grands Classiques Gauthier-Villars).

\noindent
- B. de Finetti, Funzione caratteristica di un fenomeno aleatorio,
{\it Memorie della R. Accademia Nazionale dei Lincei}, {\bf 4} No 5,  86-133 (1930).
%%[reprinted in \cite{DeFinetti OPERE06} Vol. I,
%pp. 109-158]

\noindent
- A.~Ya.~Khintchine, {\it Asymptotische Gesetze
der Wahrscheinlichkeitsrechnung}, Julius Springer, Berlin, 1933.

\noindent
-  A.~N.~Kolmogorov, {\it Grundbegriffe der
Wahrscheinlichkeitsrechnung}, Julius Springer, Berlin, 1933.

%% the appearance in 1934 of L\'evy paper on an important class of
%% stochastic processes.\

\noindent
- P.~L\'evy,
  Sur les int\'egrales dont les \'el\'ements sont des variables al\'eatoires
 ind\'ependentes,
   {\it Annali della R. Scuola Normale di Pisa} {\bf 3},  337-366 (1934)
    and  {\bf 4}, 217-218 (1935).

The first papers by Khintchine on Probability ([G~14], [G~15])
were devoted to the law of iterated logarithm for a Bernoulli sequence of random variables.

Let us, for instance, formulate in Khintchine's form the result in [G~14]: let we have an infinite series of mutually independent trials
(experiments) in any of which the probability of an appearance of
certain event $E$ is equal to $p,~0<p<1$. Suppose that the event
$E$ is realized $m(n)$ times in first $n$ trials and denote
$$
\mu(n) = m(n) - np.
$$
%It is known that with the probability arbitrary close to 1 the
%sequence $\mu(n)$ is infinitesimal with respect to $n$. 
It is
important to determine an exact upper limit for the order of
$\mu(n)$ as well as an exact sense of such an order.

{\bf Problem.} {\it To find a function $\chi(n)$ satisfying the
following conditions:

\noindent for any arbitrary small $\varepsilon > 0$ there exists a
positive integer $n_0 = n_0(\varepsilon)$ such that with
probability greater than $1 - \varepsilon$ the following
assertions hold:

$1^0.$ For all $n > n_0$
$$
\left|\frac{\mu(n)}{\chi(n)}\right| < 1 + \varepsilon.
$$

$2^0.$ There exists $n > n_0$ for which
$$
\left|\frac{\mu(n)}{\chi(n)}\right| > 1 - \varepsilon.
$$
}

{\bf Answer.} {\it The solution of the problem is given by the
formula
$$
\chi(n) = \sqrt{2 p (1-p)n \log\log\, n}.
\eqno{(1)}
$$
Any other solution is asymptotically equivalent to this one.}

The papers [G~14], [G~15] initiated the work of the Moscow probability school, and Khintchine's work in the area
asstated by Cramer (1962) \cite{Cramer_obituary}. After this work (together with the preceded to it work by E.~Borel on
the strong law of large numbers) the problem of estimation of probability of sums of
random variables occupied larger place in the investigations (see, {\it e.g.},\, Kolmogorov's
survey in (1959) \cite{MathUSSR_17-57V1}).

In 1925 Khitchine and Kolmogorov had initiated the systematic study of the convergence
of infinite series whose terms are mutually independent random variables.\footnote{See their joint
paper [G~25]:  A.Ya.  Khintchine and A.N. Kolmogorov,
Ueber Konvergenz von Reihen, deren Glieder durch den Zufall
bestimmt werden,
  {\it Mat. Sb.} {\bf 32}, 668-677 (1925).} In this paper, Khintchine proved that for countably valued random variables,
convergence of means and variances guarantees almost sure ({\it i.e.}, with probability one)
convergence of series. In order to obtain this result Khintchine used the ideas from measure
theory, namely, constructed the corresponding random variables as functions on the interval $[0, 1]$ with Lebesgue measure.
From the modern point of view it is not necessary to apply such a construction, but, in a sense,
it was one of small stones anticipated the Kolmogorov's foundation of Probability Theory.\footnote{Chapter 6 of his book
\cite{Kolmogorov_found} Kolmogorov had devoted to the results of Khintchine and himself on the applicability
of the ordinary and strong law of large numbers. In the Preface  to this book he wrote: ``I wish to express
my warm tanks to Mr. Khintchine, who has read he whole manuscript and proposed several improvements.''}

Khintchine's papers devoted to the law of the iterated logarithm and the summation of series of random terms were
followed by the papers dealing with the classical problem of summation of independent random variables.
He gave the particularly clear conditions of the applicability of the law of the large numbers
in the case of mutually independent, identically distributed summands,
which is reduced to the existence of a finite mathematical expectation ([G~44]).
Namely, he had proved that with the probability arbitrarily close to one,
$\sum\limits_{j=1}^{n} {\bf x}_j/n \rightarrow {\bf E}({\bf x}_1)$
for any sequence ${\bf x}_1, {\bf x}_2, \ldots$ of mutually independent random variables with a common distribution,
having a finite expectation.

In many of his works, special attention was paid by Khintchine to the conditions of convergence
to the Gaussian Law.
Among the first papers dealing with this topic, we have to mention the 1929 papers [G~47] and
[G~48], which initiated a
new direction of studies of so-called case of ``large deviations.''

Later, in 1935  he started  to develop the idea of the domain of attraction of the Gaussian Law (see [G~74]). This notion, which was
 introduced by P.~L\'{e}vy has been essentially extended
in the 1938 monograph by A.Ya.Khintchine ([G~92]). Let us describe the main result of the paper [G~74] again in Khintchine form,
Firsr  he gaves the following defintion: a
distribution law  $F(x)$ is said to belong to the  domain of
attraction  of
another law $G(x)$, when for the
    the sum $S_n$ of $n$ random variables, which are independent and identically distributed under the law
    $F(x)$,
there exist numbers
$\xi _n >0$, $\eta_n$ such that the distribution law
of $S_n/\xi _n -\eta_n$ tends to $G(x)$ for $n \to \infty\,. $

Since the Gauss law is surely the most important among the limiting laws in
the above  sense, it is interesting to look for
a simple criterium which allows one to recognize if
a given distribution law  belongs to the domain of attraction of
the Gauss law.

It is known that this is the case when the integral
$$ \int_{-\infty}^{+\infty}  \alpha ^2 \, dF(\alpha )$$
is finite, but there are laws for which this integral is infinite
that belong to the domain of attraction of the Gauss law.

In  [G~74] the following theorem has been proved.

{\sc Theorem}\footnote{%%
This theorem appears as {\sc Theorem} [45], pp. 192-193,
in  1938 Khintchine's book,
{\it Limit Distributions for the Sum of
Independent Random Variables}.
In the {\it Bibliographical and Historical Notes} at the end this book,
p. 114,
 Khintchine states that theorem was proved
independently and simultaneously by him,
by  P. L\'evy:   Determination g\'en\'erale des lois limites,
  {\it Compt. Rendus Acad. Sci. Paris} {\bf 203} (16), 698-700 (1936),
and  by W. Feller:
%% W. Feller,
Ueber den zentralen Grenzwertsatz  der
Wahrscheinichkeitsrechung,
  {\it Math. Z.} {\bf 40}, 521-559 (1935).
According to B.V. Gnedenko and A.N. Kolmogorov
(see p. 172 in the 1954 English translation
of their 1949 book: {\it Limit Distribution for Sums of Independent
Random Variables})
Levy's proof is referred to
  P. L\'evy:
  Propri\'et\'es asymptotiques des sommes de variables al\'eatoires
  independentes ou enchain\'ees,
   {\it J.  Math. Pures Appl.} (ser. 9) {\bf 14}, 347-402   (1935).
Gnedenko \& Kolmogorov also quote  the 1935 paper by Feller in
{\it Math. Z.} like
Khintchine, and, in addition,   a further  note by Feller, published in
{\it Math. Z.} {\bf 43}, 301-312 (1937).
}
%%%%%%%%%%%%%%%%%%%%%%%%%%       THE END OF THE FOOTNOTE %%%
{\it
The distribution $F(x)$ belongs to the domain of attraction
of the Gauss law if and only if}
$$ \lim_{x \to +\infty}
     \frac{ x^2 [1 - F(x) - F(-x)]}
   {\displaystyle \int_{-x}^{x}  \alpha ^2 \, dF(\alpha )} = 0 \,.
 \eqno{(2)}$$

The concept of the domain of attraction of the Gaussian Law was found to be highly connected with the conditions for
convergence to Gaussian of the normed sums of independent components.
In the case of identically distributed components, Khintchine found simultaneously
with P.~L\'{e}vy and W.~Feller, but independently of them,
the necessary and sufficient conditions for convergence to the normal distribution (see  [G~79]).

\section{The interaction with Paul L\'evy}

In 1928 B. de Finetti started a research regarding functions with
random increments, see
\cite{DeFinetti 29a,DeFinetti 29b,DeFinetti 29c,%%%
DeFinetti 30h,DeFinetti 31c} based  on the theory of infinitely
divisible characteristic functions, even if he did not use such
term\footnote{More details on infinitely divisible distributions can be found below in Sec. 5, and
also in our 2006 paper \cite{FMSR_M2EF}}.

%% His  results can be summarized in a number of relevant
%% theorems   (partly stated in the previous Section). As it was
%% already mentioned they are highly connected with the stochastic
% processes with stationary independent increments.
 The general case of de Finetti's problem, including also the case of
{\it infinite variance}, was investigated %% by a very different method
in 1934-35 by L\'evy \cite{Levy 34,Levy 35}, who published  two
papers in French in the Italian Journal {\it Annali della Reale Scuola
Normale di Pisa}.

At that time L\'evy (1925) \cite{Levy CP25} was quite interested in the
so-called {\it stable distributions} that are known to exhibit
infinite variance, except for the particular case of the Gaussian.
The approach by L\'evy,
 well described later in his classical 1937 book \cite{Levy ADDITION37},
is quite independent  from that of Kolmogorov, as can be
understood from footnotes in his 1934 paper \cite{Levy 34}. From
the foot\-note$^{(1)}$ we learn that the results  contained in his
paper were presented in three communications  of the Academy of
Sciences (Comptes Rendus) of 26 February, 26 March and 7 May 1934.
Then, in  the footnote $^{(6)}$, p. 339, the Author writes:
\\
{\it [Ajout\'e \'a la correction des \'epreuves] Le r\'esum\'e de
ma note du 26 f\'evrier, r\'edig\'e par M. Kolmogorov, a attir\'e
mon attention sur deux Notes de M. B. de Finetti} (see
\cite{DeFinetti 29a,DeFinetti 30h}) {\it et deux autres de M.
Kolmogorov lui-m\^eme} (see \cite{Kolmogorov 32a,Kolmogorov 32b}),
{\it  publi\'ees dans les Atti Accademia Naz. Lincei (VI ser). Ces
derni\`eres notamment contiennent la solution du probl\`eme
trait\'e dans le pr\'esent travail, dans le cas o\`u le processus
est homog\`ene et o\`u la valeur probable $\EE\{x^2\}$ est finie.
Le r\'esultat fondamental du pr\'esent M\'emoire apparait donc
comme une extension d'un r\'esultat de M. Kolmogorov.}
%, where, as noted before,
%the term {\it infinite divisible} first appeared!.

This  means that P.~L\'{e}vy was not aware  about the results on
homogeneous processes with independent increments obtained by
B.~de Finetti and  by A.~N.~Kolmogorov. The final result of L\'evy
is
%reported in Section 2 as Eq. (2.8),
known as the {\it  L\'evy
canonical representation} of the infinitely divisible
characteristic functions.

In a paper of 1937 Khintchine  [G \ref{Khintchine BMGU37a}] showed that
L\'evy' s result can be obtained also by an extension of
Kolmogorov's method: his final  result  is known as the
{\it L\'evy-Khintchine  canonical representation} of the infinitely
divisible characteristic functions. The translation from the
Russian of this fundamental paper can be found in [G~81].
The theory of infinitely distributions  was then  presented  in
German in the article [G \ref{Khintchine MatSb37}] and in Russian in
his 1938 book on {\it Limit Distributions  for Sums of Independent
Random Variables} [G \ref{Khintchine LIMITDISTR38}].

The importance of the results by A.Ya.~Khintchine has been highly
estimated by P.~ L\'evy in  the 1970 book on his life \cite{Levy LIFE70}.
At p. 105  he writes: {\it ``Tout cela constituait un ensemble important,
d'autant plus que a th\'eorie des lois ind\'efinitement divisible
a \'et\'e la base d'un autre chapitre important du calcul des
probabilites, l'arithm\'etique des lois de probabilit\'e, dont je
parlerai plus loin. Sans songer \`a minimiser le r\^ole de
Khintchine, qui m'a devanc\'e sur plusieurs points et qui san
doute aurait retrouv\'e ceux pour lesquels je l'ai devanc\'e, sans
oublier non plus les travaux ant\'erieures
 de Cauchy, de P\'olya, de B. de Finetti et de Kolmogorov, je crois  (also
pouvoir dire) que cette th\'eorie est essentiellement mon oeuvre."}

A.Ya.~Khintchine and P.~ L\'evy were working on the creation of the
theory of stable law distributions in competition\footnote{%%%
Khintchine paid a lot of attention to the results by L\'evy.
Due to him  L\'evy's theory became more accessible to the international probability community.}.
They were interested in  the works of each other. Some of their results look as a replica
to those obtained by the other side. One  example of
such  papers is a small note by A.Ya.~Khintchine in 1937 (see
[G~84]) devoted to the description of the invariant classes
of distributions. The obtained theorem is a direct generalization
of an analogous theorem by P.~L\'evy (see below). It is stipulated
the role of the invariant classes in the Probability Theory.
Another  paper is  the 1937 note by A.Ya.~Khintchine (see
[G~85]) in which examples of stable distribution laws are
constructed.
%%%
A.Ya.~Khintchine and P.~ L\'evy  had published only one joint
paper in 1937 [G \ref{KhintchineLevy CRASP37}].
It was devoted to the proof of
the following theorem on the representation formula for characteristic
functions of stable distribution laws:
%%\\
{\it The characteristic function  $\varphi(t)$
(i.e. expected value of $\e^{itX}$, $X$ being the random variable under consideration)
of the stable distribution law
is defined by
$$
\log \, \varphi(t) = -c \left( 1- i\beta \frac{t}{|t|} \, \tang \,\frac{\pi}{2}\alpha
\right) \, |t|^\alpha
\quad (c>0\,,\; 0<\alpha \le 2\,,\; |\beta| \le 1 )\,.
\eqno{(3)} $$}
%%%%
The proof of the theorem consists
of two parts. The part concerning the case $0 < \alpha < 1$ was
written by P.~ L\'evy, and that concerning the case
$1 < \alpha < 2$ was written by A.Ya.~Khintchine. It means that the paper was
written by correspondence and the authors did not meet to discuss
the results of the paper.
\newpage
\section{Infinitely divisible distributions}

In the middle of 1930s A.Ya.~Khintchine  constructed the
general theory of the limit distributions for sums of independent
random variables. In his 1937 fundamental paper [G \ref{Khintchine
MatSb37}]  he
 obtained the following main results.

 \noindent
 - He  proved the Kolmogorov's conjecture,
namely: the limit law for the sums of mutually independent random
variables, such that every summand is neglected with respect to
the sum,  have to be {\it infinitely divisible}\footnote{%%
For  details on the origin of the infinitely divisible distribution we
refer the reader to our 2006 paper \cite{FMSR_M2EF}: Mainardi, F. and  Rogosin, S., The origin of infinitely divisible
distributions: from de Finetti's problem  to L\'evy-Khintchine formula,
 {\it  Mathematical Methods for Economics and Finance}, {\bf 1} (2006) 37--55.
 [E-print {\tt http://arXiv.org/math/arXiv:0801.1910}]
 }.
%%%
It  happens also that the class of infinitely divisible
distributions is exactly the union of those distributions which
are the limits of sums of mutually independent random variables
satisfying the condition that none individual summand has an
influence on the value of the limit of sums. It is clear that the
problem has no sense without the last condition since in this case
almost all sums of such a type has no limit.

\noindent
- He discovered that each partial distribution is an infinitely divisible one.
Conversely, each infinitely divisible distribution is a partial
distribution.

\noindent
-  He showed  how one can determine a new  random
variable, having an arbitrary (but not Gaussian) stable
distribution, from a sequence of
independent identically distributed random variables (with
distribution $H(x)$) by using a simple construction.
In this manner, a model can be obtained for random variables
satisfying the stable distribution\footnote{Another model of such a type
which uses the Poisson law instead of $H(x)$ and with integrals
instead of series was proposed by P.~L\'{e}vy.}. Possible
generalizations of such constructions are discussed as well.

 The above results clarify
the fact that infinitely divisible
 distributions play an important role in the classical theory of
 summation of independent random variables.

To describe the role of the results by A.Ya.~Khintchine presented
in [G \ref{Khintchine MatSb37}] we can cite here a part of the
article by A.N. Kolmogorov  {\it Limit Theorems} in 1959 \cite{MathUSSR_17-57V1}:
%% written by A.N.~Kolmogorov for the collection of papers ``Mathematics in the
%% USSR for 40 years: 1917--1957''.
%% He wrote:
``The branch connected
with the so called {central limit theorem} on attractions of
distributions of sums of a large number  independent or weakly
dependent summands (scalar or vectorial) to the normal Gaussian
distribution was developed in some directions:

\noindent
- Sharpening the classical limit theorems on attraction by the
Gaussian distribution.

\noindent
- Investigation of this problem as a special case of the problem
on attraction by arbitrary infinitely divisible distributions. It
was understood due to appearance of the latter class of
distributions and after the proof by A. Ya. Khintchine [G
\ref{Khintchine MatSb37}] in 1937 the main theorem which states
that the limit distribution of sums of independent individually
negligible summands can only be  infinitely divisible.%%
%% {\footnote{From the Obituary of Khintichine by Gnedenko(1960/61):
%% Khintchine (Khinchin)'s construction of the general theory of
%% limit distributions for sums of independent random variables [91 =
%% this paper] belongs to the same group of ideas. The fundamental
%% proposition of the theory he developed can be formulated as
%% follows; the class of limit distributions for sums of independent
%% infinitesimal random variables coincides with the class of
%% infinitely divisible distributions. The proof of this fact, and
%% also of some other propositions of the theory of summation,
%% required the development and organization of the theory of
%% infinitely divisible laws, which were then just introduced by
%% Bruno De Finetti (1929) and A.N. Kolmogorov (1932a, 1932b) in
%% papers published in Italian on Atti della Reale Accademia dei Lincei.}}.

This second   direction was  based from the very beginning on the
idea of comparing of the process of forming the successive sum of
independent variables as the limit process with independent
increments. The latter in any case has to satisfy ``unboundedly'' divisible
distribution laws (in the Gaussian case this idea was developed by
Bachelier already in 1900).''

  We  add few words by B.V.~Gnedenko in Kintchine's obituary about the importance of the class of
 infinitely divisible distributions for the development of the
 Probability Theory in the late 1930s:
%{\it B.V.~Gnedenko in Appendix3-4.}
``The fundamental role of the infinitely divisible distributions
in the Probability Theory was established by recent investigations
of B. de Finetti, A.~N.~Kolmogorov and P.~L\'evy on homogeneous in
time stochastic processes and of A.~Ya.~Khintchine and
G.~M.~Bawly\footnote{Gregory Minkelevich Bawly  (1908-1941)
graduated at the Moscow State University in 1930, defended his PhD
thesis under guidance of A.N. Kolmogorov in 1936. His scientific
advisor had greatly esteemed his results on the limit
distributions for sums of independent random variables and cited
him in his book with Gnedenko \cite{Gnedenko-Kolmogorov
LIMITDISTR54}. G.M. Bawly lost his life  in Moscow in November
1941 at a bombing attack. According to Khintchine [G 92] %% Khintchine LIMITDISTR38
 the  name of  {\it infinitely divisible
distributions} (in  a printed version)  is  found in the 1936
article by G.M. Bawly \cite{Bawly 36}, that was recommended for
publication in the very important starting volume of the new
series of Matemati\v{c}eski Sbornik.  His second (and the last)
paper \cite{Bawly 37} was published in 1937 in Turkey and contained a
local theorem on the
 limit distribution of sums of independent random variables
generalizing the corresponding result by von Mises.} on the limit
laws for sum of independent random variables.''

The results by A.Ya.~Khintchine of [G \ref{Khintchine MatSb37}]
formed  one of the main starting points of the modern theory of
limit distributions for sum of independent random variables.

In close connection with the summation theory are the results in 1937 by
A.Ya.~Khintchine on arithmetic of the distribution laws (see [G~82]). We also mention
 the result of  Khitnchine's paper [G~83],  in which a criterion for
characteristic function is given, generalizing Bochner's type
criterion.

\section{Khintchine's book on the distribution of
the sum of independent random variables}

The 1938 monograph [G \ref{Khintchine LIMITDISTR38}] gave a
well-composed presentation of the general limit theorem for sums
of independent random variables and their application to the
classical problem on the convergence of normed sums to the normal
law on the stage known in 1938. The book was based on a course of
lectures which A.Ya.~Khintchine had delivered at Moscow State
University. This course attracted the interest
of A.A.~Bobrov, B.V.~Gnedenko, D.A.~Raikov towards the problems of
summation of random variables.

Among the results, one  finds in this book the theorem which
states that the class of all stable distributions coincides with
the class of limiting distributions for the normed sums of the
form $\left({\bf x}_1 + {\bf x}_2 + \ldots + {\bf x}_n\right)/A_n
- B_n$, where ${\bf x}_k$ are independent and identically
distributed random variables, and $A_n > 0$, $B_n$ are constants.

To the best of authors' knowledge the book by A.Ya.~Khintchine was
never translated into another language. Further development of his
results one can find  in the classical treatise by A.N.~Kolmogorov
and B.V.~Gnedenko \cite{Gnedenko-Kolmogorov LIMITDISTR54}
 on {\it Limit Distributions for Sums of Independent Random Variables}%%%%%%
\footnote{We note that the titles of both books by
Khintchine and by Gnedenko \& Kolmogorov are identical, although in
the reference list  of the book by
Gnedenko \& Kolmogorov (in both editions),  the title of Khintchine's previous book is
in some way different ({\it Limit Theorems  for Sums of
Independent Random Variables}).} that  appeared in Russian in
1949 and in English in 1954.

The complete presentation (translation) of original
techniques and results of  the 1938 book by A.Ya.~Khintchine is made in our 2010 book
\cite{RogMai_Legacy}.
It can be interesting to
the experts in Probability Theory as well as to beginners in
this branch of Mathematics since the presented results are simple, clear and
precise.

\section{Teaching Probability Theory and Analysis}

Alexander Yakovlevich Khintchine had a constant and deep interest
to the problems of teaching as in universities as in the secondary
schools. His pedagogical ideas he has presented in his textbooks,
monographs and special articles.\footnote{A number of his articles
devoted to the questions of methodics and pedagogic was published
as special chapters of the book A.Ya.Khintchine, ``Pedagogical
Articles'', ({\it Pedagogicheskie stat'i}), Prosvestchenie,
Moscow,   1963, pp. 204 (in Russian) (see [G \ref{Khintchine
Pedagogical}]).}

In 1938-1940 he headed the physical-mathematical section of
the Methodical-educational Soviet at the Ministry of Education of
the Russian Federation. When the Academy of Pedagogical Sciences
of the Russian Federation was founded, he became an academician of
this Academy. He was very active as a member of the editorial board of
the multi-volume ``Encyclopedy of Elementary Mathematics'', some
volumes of which appeared in the late 1950s
(see e.g. [G \ref{KhintchineYaglom CHILDREN59}]).

In the article ``Alexander Yakovlevich Khintchine''
B.V.~Gnedenko says (see [G \ref{Khintchine Pedagogical}, pp. 180-196]):
``Everybody who could hear either a mathematical course lectured
by
 or his scientific report remembered how exceptionally careful
 was formulations of the problems, how exact was the description
 of the place of these problems in the corresponding branch of science.
 The reader could forget that his ideas was alien to him at the beginning.
 He/she started to feel the importance and significance of the discussed
 ideas, and to have satisfaction from the fact that he/she could
 achieve a new level of knowledge and a real understanding
 of the complicate  notions and methods. The
  success of A.Ya.~Khintchine as a Teacher and a Creator of new
  Research Direction in the Probability Theory, in the Number Theory,
  and in the Real Analysis is, in particular, based on such his way of thinking.''

  From this point of view, it would be natural to include
  into this section of our paper  certain sections and chapters
  from his monographs in which one can see the hand of a real
  Master. But in this case this section would be overloaded.

We restrict ourselves with the discussion of ideas by A.Ya.~Khintchine, which concern  teaching in
universities and in the secondary schools.

 \vspace{3mm}
 \noindent {\bf\large 7.1 Pedagogical credo}
 \vspace{3mm}

In our previous sections, one can see
all features of texts written by A.Ya.~Khintchine. Many of his
books run into several editions. We have to mention here few words
by B.V.~Gnedenko from the Introduction to the 4-th edition in 1978 of the
book in Russian {\it Continued Fractions} [G \ref{Khintchine
CONTFRACT35-49-63}]:
``Before starting to  write this Introduction
I have read the book once more. I feel the great pleasure of the
contact with the big master again. He does not only perfectly know
the material but also can present it in such a way that a reader
find himself under the influence of the author's individuality.
His book is not consist simply of declarations but it contains
complete proofs of all statements. Anyway (since the book is
addressed mainly to the beginners), all the proofs help the reader
to understand the real course of reasoning, as well as their
necessity.''
Such sentence means that in all his books
A.Ya.~Khintchine is not only a scientist, but also, and mostly,
 a teacher.

Another example of the same type we can meet by reading the small
booklet ``{\it Three Pearls of Number Theory}''  which has been appeared  in
several Russian editions (see [G \ref{Khintchine NUMTHEORY47}]). It is
written in the form of a letter to a former student  of Moscow
University (a soldier of the second World War), who asked the author to send some
mathematical pearls to him in hospital. Even in this case
A.Ya.~Khintchine provided the reader with three very interesting
results on number theory which were proved by rather young
mathematicians.
The first pearl is concerning with the hypothesis on
arithmetic progression which was solved by the young Dutch mathematician B.L.~van der Warden in 1928.
The second pearl deals with the  hypothesis of Landau-Shnirelman on the thickness
of sums of sequences which was solved by the young American mathematician  G.~Mann in 1942.
The third pearl is the elementary proof of  Waring's  problem
in Logic given by the young Russian mathematician
Yu.V.Linnik in 1942.
Khintchine told the story of these
proofs by describing in detail (but on quite elementary level)
all their steps. Thus he tried to involve the reader into the
process of thinking saying him: you see, these people were of the
same age as you and they formulated and successfully attacked very
serious mathematical questions by using deep but simple
considerations.

Alexander Yakovlevich Khintchine was one of the best lecturers of
Moscow State University. During several decades he delivered a
course of Mathematical Analysis in universities and pedagogical
institutes. His students considered these lectures as the best
they have ever attended during  their students' years. It was no
wonder in it. A.Ya. Khintchine tried to carry out for listeners
not only the formal mathematical machinery, not only analytical
technique, but the essence and ``soul'' of Analysis, its basic
ideas. His pedagogical credo was: ``better not too much but in the
best form''.
%%%%%%%%%%%%
\newpage
%%%%%%%%%%%
\vspace{3mm} {\bf\large 7.2 Mathematics in the secondary school}
\vspace{3mm}

Alexander Yakovlevich Khintchine was very enthusiastic in
strengthening the content of course of Mathematics in
secondary schools, as well as in developing new teaching methods
and methodical receipts.

He had hardly working on the problems of education at the
secondary school. A collection of his papers on these problems was
published in 1963 in Moscow [G \ref{Khintchine Pedagogical}]. In
his materials prepared to publication one can find different ideas
which concern the questions of mathematical education. For example,
he said (see [G \ref{Khintchine Pedagogical}, p. 11]): ``there are
several types of repetition or revising of the material which
pupils have to make, namely, 1) repetition before the start of
academic year, 2) revising following its beginning, 3) review of
lessons learned
 during the latter, 4) revising of the material of
certain theme connected with the control of knowledge, 5) annual
repetition, 6) revising the material at the preparation to exams.
This is awful. Whether it is possible to teach in such a way that pupils
will not forget the material in order to avoid such infinite chain
of repetitions?'' We see here the great concern of A.Ya.
Khintchine on the important educational problem, saying that
repetition is not a way out. The teacher have to find a better
method making mathematical knowledge of pupil real and useful
instrument in their life. The main requirements of A.Ya.
Khintchine to teaching of mathematics in the secondary school
are: a) an account of the age-specific features of the pupils can
lead to the necessity of simplified presentation of ideas and
notions of the science. But such  simplification should not
falsify and even distort the scientific content  of these notions;
b) replacement of the rigorous and exact definitions and proofs by
fuzzy presentations having no exact sense cannot facilitate the
understanding of the subject. Khintchine said that fuzzy thinking
is not  simpler than precise thinking.

 \vspace{3mm} {\bf 7.2.1
General ideas on mathematical education in secondary schools}
\vspace{3mm}

A.Ya. Khintchine thought that the main goal in the methodology
of mathematics as well as in the educational process as a whole is to
awake the creative ideas, to develop the technique mostly suitable
for it. For him the success of pedagogical process is not in the
collection of good marks, but in the depth of understanding.
%%of a teacher by his pupils.
Pupils have to use correct and rigorous
logical considerations, to see the gaps in arguments. The goal of
teachers is to show the ways to solve independently non-standard
problems, to sort out their known proofs.

He said ((see [G \ref{Khintchine Pedagogical}, pp. 29-30]): ``Two
principles I have to lay into the base of the solution of the
question at what level one or another mathematical notion should
be studied in the secondary school taking into account the age
features of pupils and the modern scientific understanding of this
notion. These principles are:

1. If age-specific features of the pupils do not allow to give to
certain notion its real scientific interpretation then the
conception of this notion can be simplified. It means that the
school should not develop the notion to the level accepted in the
modern mathematical science. It can be stopped at one of previous
levels of its development. But in any case the teacher should not
distort the scientific meaning of the notion to give it feature
which contradict this meaning. The school should not develop any
notion in the direction deviating from the way of its scientific
development.

2. If one changes of sharp and exact definitions, formulations and
reasoning by fuzzy ones, having no exact sense, leading at the
sequential application to the logical contradiction, then it
cannot simplify the real understanding. Fuzzy thinking cannot be
simpler than the sharp one.

We suppose also that the usual theoretical scheme of the course of
mathematics in the secondary school contains a lot of archaic
notion which stay out of the the main stream of mathematical
development. In the most cases creating specially for the school
 notions which are not used in the science has no methodical sense
 and bring the only damage to real mathematical development of
 pupils.''

 \vspace{3mm}
 \noindent
 {\bf 7.2.2 Basic mathematical
notions in the secondary school}
\vspace{3mm}

Khintchine  considered the following
mathematical notions: number, limit and function as the basic ones.
Starting the description of numbers the teacher has to formulate
and clarify the definition of this notion. He/she has to create
among his pupils the real understanding of a number as the subject
of arithmetic operations. Of course it should be done carefully,
step-by-step starting with natural numbers and integers, continued
by fractions (rational numbers). It is important to make the whole
picture of developing the notion of a number. Special attention
has to be paid to introduction of irrational numbers since it
allows to do the steps creating the continuous line
measure, description of continuous measuring of many physical
magnitudes. Finally, on this base pupils could understand the
notion of continuity, of limit and continuous function (see [G
\ref{Khintchine Pedagogical}, pp. 44-49]). Complex numbers have to
be also an element of mathematical education in the secondary
school. Teacher has to show natural is to extend real numbers in
this way, how many connections has such a notion with different
mathematical problems of algebra and geometry. It can be done by
using different geometric illustrations.

Since the notion of limit was changing in the course of its
developing one has to take it into account. This notion have to be
developed in the pupils' mind too. It is not successful to start
already with the formal modern definition. Every teacher in
mathematics knows that such formalism kills many important
questions dealing with this notion: why it is appeared, how it is
constructed, how it is connected with another basic notions of
mathematics.

Almost all methodists of Khintchine's time (as well as of our time)
had supposed that the notion of function has to be a central
pivot of the whole course of mathematics. Basing on it  the main notion of arithmetic, algebra, geometry and
trigonometry have to be created.
It is really correct but can lead to certain
overestimation and even abuse. This notion has to be applied very
carefully accounting the age of pupils, the nature of the
considered problems etc. Their study must not consist only of
drawing graphs and considering formulas. The teacher has to show
the real contents of this notion, show how can can use it for the
solution of mathematical problems. From the other side, there is no need
 to make unnecessary generalizations and extensions, such
as, e.g., introduction of multi-valued functions.

The question on the most suitable form of introduction of one or
another mathematical notion in the course of mathematics in the
secondary school is one of the main aims of the mathematical
methodology. Before starting a discussion concerning certain
concrete notion one has to understand the general principle
characteristic for the corresponding branch of mathematics. First
of all we have to answer the question what is a mathematical
definition, what is its role for the secondary school, how and in
which way one has to replace the definition in the case when its
complete introduction is logically impossible or methodically
unreasonable.

In the paper ``On mathematical notions in the secondary school''
(see [G \ref{Khintchine Pedagogical}, pp. 85-105])
A.Ya.~Khintchine tried to clarify under which conditions these notions should be
posed. He answering the questions on the most suitable introduction of
concrete mathematical notions. The aim of that paper was to create
the real background for further discussions of these questions in
order to make these discussions really productive.

\vspace{3mm}
\noindent {\bf 7.2.3 Formalism in mathematical education in the
secondary school}
\vspace{3mm}

Khintchine  considered the formalism of mathematical
knowledge and abilities as the most heavy disadvantages of
mathematical education in the secondary school. Those pupils which
get only formal part of mathematical methods became weak in front
of problems arising in the real life. They could not pose the
questions in mathematical sense and from mathematical point of
view, and could not solve the problems of such a kind.

These pupils are also in a weak positions in the universities. The
courses of mathematics in the universities are even more free from the
formal approaches, their study need to understand the corresponding
ideas rather than learn by heart the formal conceptions and
mathematical notions.

One more dangerous consequence of the formalism in the study of
mathematics in the secondary school is that the formal knowledge of
mathematics is useless for creating of scientific world-view of
the pupils.

How to fight against formalism in the concrete situation, in the
framework of the concrete programs of mathematical courses? These
questions are discussed by Khintchine on the base of some
examples taken from programs of mathematics and from
his teaching practice (see his article ``Formalism in teaching of
mathematics in the secondary school'' in [G \ref{Khintchine
Pedagogical}, pp. 106-127]).

\vspace{3mm}
\noindent {\bf 7.2.4 Educational effects of mathematical
classes}
\vspace{3mm}

The article with such title (see  [G \ref{Khintchine Pedagogical},
pp. 128-160]) was published already after the death of
A.Ya.~Khintchine though it was presented by him at one of the
scientific meetings of the Mathematical Office of the Research
Institute on school educations of the Ministry of Education of
Russian Federation. Some of his ideas, presented below, seems
rather modern. Let us present few of them.

 The objects of Mathematics (as a Science and as
a Discipline studied in the secondary school) are quantitative
relations and spatial forms of real things, but not these things
themselves. In a sense it decreases an educational effect of the
course of mathematics in the school.

The most known methods to reach such an effect are
the use of specific mathematical logical rigor to the creating
the general logic and so called culture of thinking. From the
other side, a teacher can equip certain mathematical problems by
the concrete content and it gives possibility to extend the mental
outlook of pupils, and to increase their general cultural level.
But it is not all possible.

First of all, we have to speak about creating the culture of
thinking. The teacher has to train pupils to be correctly
thinking. He should always use the complete argumentation. Then,
he has to show on the concrete examples how wrong is the way of
illegal generalizations. The teacher has to fight against
inconsistent analogies. It is important always to follow in the
proofs the complete argumentation. The creating
classifications has to be also complete and consistent.

Another moment, which is important to develop, is the concrete and
logical style of thinking. In any case it is very useful in any
area of the future life of pupils.

\vspace{3mm}
\noindent {\bf 7.2.5 On so called ``Problems on Reasoning'' in
the Course of Arithmetics}
\vspace{3mm}

Even in the modern school we can meet the idea to start teaching
pupils from the early age on the base of so call problems on
reasoning. A.Ya.~Khintchine discussed this idea basing on the
material of the course of mathematics (arithmetic) for 5th year of
the school (see  [G \ref{Khintchine Pedagogical}, pp.
161-172]).\footnote{This article was never published and was found
by B.V.~Gnedenko in the archive by A.Ya.~Khintchine.}

He said: ``If we ask even very good teacher of mathematics how
many pupils (of e.g. 5th year of the school) can solve the
problems which consist not in simple calculations but need to find
a special way of solving, then we have got not very impressive
statistics. Of course quite a lot of pupils can do it if before
they solve a number of analogous problems. Thus the aim to develop
the quick-wittedness of pupils on the base of hard and
non-standard problems can not be reached even by the very good
teachers. If we consider the concrete problems proposed for pupils,
we can see how difficult for them is to find best way for solving,
and even more difficult for teachers is to find a methodic of
creating the corresponding skills. Besides, these problems
appeared to be solving by using much more simple technique few
years later. We can not understand for what reason pupils have to
invent an unusual and nonstandard for his/her knowledge methods if
he/she will reach the solution by standard technique at the
corresponding age. Besides, even solution of quite a simple
problem which is obtained correctly and very fast can be made
pupils more happy than heavy ``creative'' thinking leading to the
solution which is impossible to understand completely.''

\vspace{3mm}
\noindent {\bf\large 7.3 Teaching mathematical analysis}
\vspace{3mm}

This part is based on the article by
A.I.~Markushevich ``A.Ya.~Khintchine as a teacher of mathematical
analysis'' included into the book A.Ya.  Khintchine,
 ``{\it Pedagogical Papers}'' ([G \ref{Khintchine Pedagogical}, pp.
 173--179]).

 A.Ya.~Khintchine had  outstanding pedagogical skills. His
 courses were short but free of unessential details. He has always tried
 to explain the importance and mathematical essence of considered
 notions and to motivate the formulation of the problem. Thus he
 prepared the audience to discovering of a new research material.
This was a necessary condition for them to follow the lecturer's
guidance and to go through complicated mathematical constructions.
From time to time he made some departure from the main course and
explain to the students how to achieve one or another pedagogical
aim. He taught them to be teachers.

He was never afraid to attack very complicated pedagogical
problems. The characteristic example of it is the well-known 1943-1948  book
``Eight lectures on Mathematical Analysis'' ([G \ref{Khintchine
MATHANAL43}]). This book is addressed to those who adopting the
machinery of the mathematical analysis tries to understand its
main ideas and logic. These 8 lectures (chapters) in the book are
the following: I. Continuum, II. Limits, III. Functions, IV.
Series, V. Derivative, VI. Integral, VII. Series representation of
functions, VIII. Differential equation. Each of this subject is
discussed on the historical retrospective, but mainly the
presentation is made in such a way which, that could form the modern
point of view of the reader.\footnote{It is amazing that after
more than 50 years from the first edition of the book it remains
to be modern and can teach young people with completely new
mentality and taught by another mathematical program in the
secondary school.}

Another brilliant example of realization of pedagogical credo by
A.Ya.~Khinchine is his 1953-1957 ``Short course of mathematical analysis''
([G \ref{Khintchine ANALYSIS53-57}]). ``Eight lectures'' are more
or less free in style and contain a lot of ``considerations'' and
explanations, but the ``Short course'' is mainly rigorous. It has
presented all necessary facts and their proofs. Another feature of
the book is that it step by step develops very vague and intuitive
mathematical notions brought by students from the secondary
school. One of the most characteristic examples is the
introduction (on a new level) of the notion of the limit. It is
done in the same way as we have described before, namely, several
stages from the school-type notion to the rigorous should be
made by the reader under the guidance of the author.

He discussed many special problems of mathematical education.
Among them is, for instance, the question of organization of the
independent and unassisted work of students (how modern for
Russian high school it sounds!). A.Ya. Khintchine said: ``When we
are talking about independent work of students we mean first of
all their consulting by teachers. In reality it means for any
difficulty (even small one which can be clarified in 10 minutes of
thinking) the student appeals to the teacher and got a final
answer or exact citing in the corresponding of the book without
any afford from the student's side. Is it a stimulation of
independent work of students?''

Special attention has to be paid from his point of view to
education and re-qualification of teachers of the secondary school.
They should get really scientific knowledge in all themes studied
in the secondary school and even more. It should not be a
collection of certain receipts how to show pupils a solution of
one or another problem. Such collections are widely spread (up to
now!)

A.Ya. Khintchine was one of the supporters of the idea to
introduce elements of analysis of infinitesimals already in the
course of the secondary school. Let as recall few of his ideas
concerning this question. ``Analysis of infinitesimals (or
Calculus) is one of the most import discovering of the mankind. It
has numerous applications. Calculus is important for formation of
the scientific point of view of our pupils.''

\vspace{3mm}
\noindent {\bf\large 7.4 Teaching probability theory}

\vspace{3mm} The contribution of A.Ya.~Khintchine in the theory of
functions and the number theory is great. Anyway, the most of his
scientific life was connected with Probability Theory. In 1920s
A.Ya.~Khintchine, A.N.~Kolmogorov, E.E.~Slutskii and P.~L\'evy have
discovered the tight connection between Probability Theory and the
mathematical branch which studied sets and general notion of the
function. Very close to the understanding of this connection had
came a little bit earlier E.~Borel.

The stochastic processes were discovered due to fundamental works
by A.N.~Kolmogorov and A.Ya.~Khintchine. In a sense this theory
was developing the ideas of A.A.~Markov to study dependent random
variables (called later Markov chains).

Systematic using of the methods of the set theory and the theory
of functions of real variables in probabilistic models,
construction of the basis of the theory of stochastic processes,
extended developing of the summation theory for independent random
variables, as well as introduction of the new approach to the
problems of statistical physics, all these thing constitutes the
important impact of A.Ya.~Khintchine to the creation of modern
Probability Theory.

Starting from the problems connected with the number theory (The
Law of Iterated Logarithm) and the theory of functions
(convergence of series of random summands) he extended his
interest to more and more deeper problems of the developing
theory. Moreover, he attracted many young Moscow mathematicians to
study these problems and thus create highly developed Moscow
school of Probability. Such behavior was really characteristic
feature of his pedagogical talent. He was really great Teacher. It
is remarkable that his only textbook on Probability (joint with
his student B.V.~Gnedenko) he entitled `Elementary Introduction to
the Probability Theory'' ([G \ref{GnedenkoKhintchine PROB46}]),
but which is neither elementary nor simple.

We have to mention also, that all his books appeared due to courses
of lectures which he was delivered in universities. Thus the first
monograph ([G \ref{Khintchine PROBABILITYT27-32}]) on Basic Laws
of Probability appeared in 1927-1932 due to his intention to understand the
nature of probabilistic approach and to show it to the young
people on the most simple examples. The second monograph deals with
his explanation of the just completed by A.N.~Kolmogorov and
I.G.~Petrovsky theory of Markov processes ([G \ref{Khintchine
ASYMPROB33-36-48}]). At last the third monograph ([G
\ref{Khintchine LIMITDISTR38}]) in 1938 is a real course of lectures. 
This course  attracted
to the problems of summation A.A.~Bobrov, B.V.Gnedenko and
D.A.Raikov.

This list is really incomplete. In any of his books (even a very small
booklet) it can be recognized his intention to widely describe the
corresponding ideas and to involve people in solving the appeared
problems. Such approach we can meet in his later works on queueing
theory (or how he was called it ``the theory of mass service''),
information theory, foundation of statistical physics. It is very
interesting to read paper by A.Ya.~Khintchine in which he discuss
the main ideas of probability with school boys and girls (see
article in ``Children Encyclopedia'' [G \ref{KhintchineYaglom
CHILDREN59}] and brilliant booklet [G \ref{Khintchine CHANCE34}]).

He  supposed that elements of probability have to be taught
already in the secondary school. His main arguments were: 1) many
people finished the secondary school will encounter  certain statistical data; 
the study of certain notions of the
probability theory is a good base for it; 
2) solving concrete problems of the probability theory by using formulas of
combinatorics, the  pupils may find an interest in these formulas and try to
understand on a better level the meaning of them.

\section*{Acknowledgement} This work by S.R. is supported by the Belarusian
Fund for Fundamental Scientific Research through grant F17MS-002, by ISA 
(Institute of Advanced Studies of the Bologna University)
and by Erasmus Programme which helps authors to provide a joint
research at the Bologna University. T
he work of F.M. has been carried out in the framework of
the National Group of Mathematical Physics (GNFM--INdAM) activities.

\newpage

\begin{center}
%\vvs
\medskip
 \def\date#1{\gdef\@date{#1}} \def\@date{\today}
%% {\bf Ver Final, \@date }
%%  {PRELIMINARY VERSION \@date}
%%%%%%%%%%%%%%%%%%%%%%

\end{center}

%\newpage

\centerline{{\bf BIBLIOGRAPHY on A. Ya. KHINTCHINE (A. I.
KHINCHIN)}} \vspace{0.25truecm}
 \centerline{{\bf from his Obituary by B.V. Gnedenko, University of Kiev}}
\vspace{0.25truecm} \centerline{{\bf  The 4-th Berkeley Symposium
 on Mathematical Statistics and Probability}}
\vspace{0.25truecm} \centerline{ {\bf held at the Statistical
Laboratory, University of California}}

\vspace{0.25truecm} \centerline{{\bf Berkeley, June 20 - July 30,
1960}}

\vspace{0.5truecm} \centerline{{\bf Proceedings, Vol. II,
University of California Press,}} \vspace{0.25truecm}
\centerline{{\bf Berkeley and Los Angeles, 1961, pp. 10-15.}}

%\begin{center}
% {\rm\bf\ Список публикаций соискателя \par}
%\end{center}

\newcounter{N}
\begin{list}{[G~\arabic{N}]}{
\usecounter{N}}

%\item\label{RogVaitKiev} Вайтехович,~Т.С. Ограниченность оператора
%сужения аналитической функции / Т.С.~Вайтехович, С.В.~Рогозин //
%Сборник трудов института математики НАН Украины. -- 2006. -- Т.~3,
%№~4. -- С.~1--10.

\vspace{0.5truecm}

%%%%%%  REFERENCES

%\begin{thebibliography}{99}

%%%

%% khintchine_2.ref

%% References on KHINTCHINE (KHINCHIN) present in  fmbiblio.tex
%% and in  the bok by Gnedenko-Kovalenko and in the Obituary by Gnedenko
%%% Rivised version of this list is presented in the book
%%% A.Ya.Khintchin, Pedagogical articles, ed. by B.V.Gnedenko,
%%% Moscow: Academy of Pedagogical Sciences of Russian Federation, 1963, 204 p.
%% In vista della visita di Rogozin, portare a Pinarella
%% i libri/articoli  di KHINTCHINE (vedi folder da PIVA) e di
%% GNEDENKO-KOLMOGOROV

\item\label{Khintchine CRASP16}  %% 1 %%
 A.Ya.  Khintchine,
Sur une extension de l'int\'egrale de M. Denjoy,
  {\it C.R. Acad. Sci.  Paris} {\bf 162}, 374-376 (1916).

\item\label{Khintchine CRASP17}   %% 2 %%
 A.Ya.  Khintchine,
Sur la d\'erviation asymptotique,
  {\it C.R. Acad. Sci.  Paris} {\bf 164}, 142-145 (1917).

\item\label{Khintchine MatSb18}     %% 3 %%
 A.Ya.  Khintchine,
On the process of integration of Denjoy, {\it Mat. Sb.} {\bf 30},
548-557 (1918) [in Russian]. %!%

\item\label{Khintchine IIVPI21}     %% 4 %%
 A.Ya.  Khintchine,
Sur la th\'eorie de  l'int\'egrale de M. Denjoy,
 {\it Izv. Ivanovo-Voznesenskogo Polytech. Inst.}
{\bf 3}, 49-51 (1921).

\item\label{Khintchine IIVPI22a}     %% 5 %%
 A.Ya.  Khintchine,
On a property of continued fractions and its arithmetic
applications,
  {\it Izv. Ivanovo-Voznesenskogo Polytech. Inst.}
{\bf 5}, 27-41 (1922) [in Russian].

\item\label{Khintchine IIVPI22b}     %% 6 %%
 A.Ya.  Khintchine,
On the question of the representation of a number as the sum of
two primes,
  {\it Izv. Ivanovo-Voznesenskogo Polytech. Inst.}
{\bf 5}, 42-48 (1922) [in Russian].

\item\label{Khintchine IIVPI22c}    %% 7 %%
 A.Ya.  Khintchine,
A new proof of the fundamental theorem in the metric theory of
sets,
  {\it Izv. Ivanovo-Voznesenskogo Polytech. Inst.}
{\bf 6}, 39-41 (1922) [in Russian].

\item\label{Khintchine FM23a}    %% 8 %%
 A.Ya.  Khintchine,
Sur les suites des fonctions analytiques born\'ees dans leur
ensemble,
  {\it Fundamenta Mathematicae} {\bf 4}, 72-75 (1923).

\item\label{Khintchine FM23b}       %% 9 %%
 A.Ya.  Khintchine,
Das Steitigkeitsaxiom des Linearkontinuums als Induktionsprinzip
betrachtet,
  {\it Fundamenta Mathematicae} {\bf 4}, 164-166 (1923).

\item\label{Khintchine MathZ23a}    %% 10 %%
 A.Ya.  Khintchine,
Ueber dyadische Br\"uche,
  {\it Math. Z.} {\bf 18}, 109-116 (1923).

\item\label{Khintchine MathZ23a}    %% 11 %%
 A.Ya.  Khintchine,
Ein Zatz \"uber Kettenbr\"uche mit aritmetischen Anwendungen,
  {\it Math. Z.} {\bf 18}, 289-306 (1923).

\item\label{Khintchine MatSb24a}    %% 12 %%
 A.Ya.  Khintchine,
On sequences of analytic functions,
  {\it Mat. Sb.} {\bf 31}, 147-151 (1924) [in Russian].

\item\label{Khintchine MatSb24b}      %% 13 %%
 A.Ya.  Khintchine,
Investigations on the sructure of measurable functions, ch.1, %!%
  {\it Mat. Sb.} {\bf 31}, 265-285 (1924) [in Russian].

\item\label{Khintchine CRASP24}         %% 14 %%
 A.Ya.  Khintchine,
Sur une th\'eor\`eme g\'en\'eral relatif aux probabilit\'es
d\'enombrales,
  {\it C.R. Acad. Sci.  Paris} {\bf 178}, 617-619 (1924).

\item\label{Khintchine FM24}         %% 15 %%
 A.Ya.  Khintchine,
Ueber einen Satz der Wahrscheinlichkeitsrechnung,
  {\it Fundamenta Mathematicae} {\bf 6}, 9--20 (1924).

\item\label{Khintchine MathAnn24}    %% 16 %%
A.Ya.  Khintchine, Einige S\"atze \"uber Kettenbr\"uche, mit
Anwendungen auf die Diophantischen
 Approximationen,
  {\it Math. Annalen} {\bf 92}, 115-125 (1924).

\item\label{Khintchine IIVPI25}         %% 17 %%
 A.Ya.  Khintchine,
On a question in the theory of Diophantine approximations {\it
Izv. Ivanovo-Voznesenskogo Polytech. Inst.} {\bf 8}, vyp. 2, 32-37
(1925) [in Russian]. %!%

\item\label{Khintchine MathZ25}        %% 18 %%
 A.Ya.  Khintchine,
 Zwei Bemerkungen zu einer Arbeit des Herrn Perron,
  {\it Math. Z.} {\bf 22}, 274-284 (1925).

\item\label{Khintchine MatSb25a}      %% 19 %%
 A.Ya.  Khintchine,
Investigations on the sructure of measurable functions, ch. 2,
 {\it Mat. Sb.} {\bf 32}, 377-433 (1925) [in Russian]. %!%

\item\label{Khintchine MatSb25b}       %% 20 %%
 A.Ya.  Khintchine,
Ueber die angen\"aherte Aufl\"osung linearer Gleichungen in ganzen
Zahlen, {\it Mat. Sb.} {\bf 32}, 203-219 (1925).

\item\label{Khintchine MatSb25c}               %% 21 %%
 A.Ya.  Khintchine,
Zur Theorie der diophantischen Approximationen,
 {\it Mat. Sb.} {\bf 32}, 277-288 (1925).

\item\label{Khintchine MatSb25d}           %% 22 %%
 A.Ya.  Khintchine,
Bemerkungen zur metrischen Theorie der Kettenbr\"uche,
  {\it Mat. Sb.} {\bf 32}, 326-329 (1925).

\item\label{Khintchine MathZ25}                 %% 23 %%
 A.Ya.  Khintchine,
Bemerkung zu meiner Abhandlung `Ein Satz \"{u}ber
  Kettenbr\"uche mit arithmetischen Anwendungen',
  {\it Math. Z.} {\bf 22}, 316 (1925). %!%

\item\label{Khintchine MatSb25e}            %% 24 %%
 A.Ya.  Khintchine,
On the Petersburg game, {\it Mat. Sb.} {\bf 32}, 330-341 (1925).

\item\label{KhintchineKolmogorov MatSb25}        %% 25 %%
 A.Ya.  Khintchine and A.N. Kolmogorov,
Ueber Konvergenz von Reihen, deren Glieder durch den Zufall
bestimmt werden,
  {\it Mat. Sb.} {\bf 32}, 668-677 (1925).

\item\label{Khintchine MatSb25e}                %% 26 %%
 A.Ya.  Khintchine,
Ueber die Anwendbarkeitsgrenzen des Tschebycheffschen Satzes in
der Wahrscheinlichkeitsrechung,
  {\it Mat. Sb.} {\bf 32}, 678-687 (1925).

\item\label{Khintchine MathZ26}               %% 27 %%
 A.Ya.  Khintchine,
 Zur metrischen Theorie der diophantischen Approximationen,
{\it Math. Z.} {\bf 24}, 706-714 (1926).

\item\label{Khintchine MathAnn26}          %% 28 %%
 A.Ya.  Khintchine,
 Ueber das Gesetz der grossen Zahlen,
  {\it Math. Annalen} {\bf 96}, 152-158 (1926).

\item\label{Khintchine Palermo26}             %% 29 %%
 A.Ya.  Khintchine,
 Ueber eine Klasse linearer diophantischer Approximationen,
  {\it Rendiconti del Circolo Matematico di Palermo}
{\bf 50}, 170-195 (1926).

\item\label{Khintchine VestKommAkad26}       %% 30 %%
 A.Ya.  Khintchine,
Ideas of intuitionism and the struggle for the subject in
contemporary mathematics,
  {\it Vestnik Kommunist. Akad.} {\bf 16}, 184-192 (1926) [in Russian].

\item\label{Khintchine FM27}                     %% 31 %%
 A.Ya.  Khintchine,
Recherches sur la structure des fonctions mesurables,
  {\it Fundamenta Mathematicae} {\bf 9}, 212--279 (1927).

\item\label{Khintchine MatSb27}                  %% 32 %%
 A.Ya.  Khintchine,
 Ueber  diophantintische Approximationen h\"oheren Grades,
  {\it Mat. Sb.} {\bf 34}, 109-112 (1927).

\item\label{Khintchine RussianCongr27}         %% 33 %%
 A.Ya.  Khintchine,
 Diophantintine Approximations,
 In: {\it Trudy All-Russian Math. Congress}, 131-137 (1927).

\item\label{Khintchine FERMAT27-32}            %% 34 %%
 A.Ya.  Khintchine,
 {\it Great Theorem of Fermat},
Gosizd, Moscow, 1927; 2-nd Edition GTTI, Moscow, 1932, 52 p. [in
Russian].

\item\label{Khintchine PROBABILITYT27-32}      %% 35 %%
 A.Ya.  Khintchine,
 {\it Fundamental Laws of Probability},
 Moscow State University, Moscow, 1927; 2-nd Edition GTTI, Moscow,
  1932, 82 p. [in Russian] %% [X]

\item\label{Khintchine CRASP28}                %% 36 %%
 A.Ya.  Khintchine,
Sur la lois forte des grands nombres,
  {\it C.R. Acad. Sci.  Paris} {\bf 186}, 285-287 (1928).

\item\label{Khintchine BELGIUM28}             %% 37 %%
 A.Ya.  Khintchine,
Objection \'{a} une note de M.M. Berzin et Errera,
{\it  Acad. Roy. Belg. Bull. Cl. Sci.,} [Ser 5] {\bf 14}, 223-224 (1928).

\item\label{Khintchine VestStat28}         %% 38 %%
 A.Ya.  Khintchine,
 Strong law of large numbers and its significance in mathematical statistics,
{\it Vestnik Stat.} {\bf 29}, 123-128 (1928) [in Russian].

\item\label{Khintchine IzvMoscow28}      %% 39 %%
 A.Ya.  Khintchine,
Begr\"undung der Normalkorrelation nach der Lindeberg\-schen Methode,
{\it Izv. Ass. Research Inst. Moscow Univ.}  {\bf 1}, 37-45 (1928).

\item\label{Khintchine MatSb28a}         %% 40 %%
 A.Ya.  Khintchine,
 Ueber  die Stabilit\"at zweidimensionaler Verteil\-ungs\-gesetze,
{\it Mat. Sb.} {\bf 35}, 19-23 (1928).

\item\label{Khintchine MatSb28b}      %% 41 %%
 A.Ya.  Khintchine,
 Ueber die angen\"aherte Auflo\"sung linearer Gleichungen in ganzen  Zahlen,
  {\it Mat. Sb.} {\bf 35}, 31-33 (1928).

\item\label{Khintchine MatSb28c}        %% 42 %%
 A.Ya.  Khintchine,
Theory of numbers: outline of its developmens in 1917-1927,
  {\it Mat. Sb.} {\bf 35} (suppl. issue), 1-4 (1928).

\item\label{Khintchine VestKommAkad29}     %% 43 %%
 A.Ya.  Khintchine,
The role and the character of induction in mathematics,
  {\it Vestnik Kommunist. Akad.} {\bf 1}, 5--7 (1929) [in
  Russian]. %!%

\item\label{Khintchine CRASP29a}        %% 44 %%
 A.Ya.  Khintchine,
Sur la lois  des grands nombres,
  {\it C.R. Acad. Sci.  Paris} {\bf 188}, 477-479 (1929).

\item\label{Khintchine CRASP29b}        %% 45 %%
 A.Ya.  Khintchine,
Sur une gen\'eralisation des quelques formules classiques,
  {\it C.R. Acad. Sci.  Paris} {\bf 188}, 532-534 (1929).

\item\label{Khintchine MathZ29}       %% 46 %%
 A.Ya.  Khintchine,
 \"Uber ein Problem der Wahrscheinlichkeitsrechnung,
  {\it Math. Z.} {\bf 29}, 746-752 (1929).

\item\label{Khintchine MathAnn29a}         %% 47 %%
 A.Ya.  Khintchine,
 Ueber die positiven und negativen Abweichungen des arithmetischen Mittels,
  {\it Math. Annalen} {\bf 101}, 381-385 (1929).

\item\label{Khintchine MathAnn29b}         %% 48 %%
 A.Ya.  Khintchine,
 Ueber einen neuen Grenz\-wertsatz der Wahrsche\-inlichkeits\-rechung,
  {\it Math. Annalen} {\bf 101}, 745-752 (1929).

\item\label{Khintchine MatSb29}             %% 49 %%
 A.Ya.  Khintchine,
 Ueber Anwendbarkeitskriterien f\"ur des Gesetz der grossn Zahlen,
  {\it Mat. Sb.} {\bf 36}, 78-80 (1929).

\item\label{Khintchine UspFiz29}        %% 50 %%
 A.Ya.  Khintchine,
The teachings of von Mises on probability and the principles of
statistical physics,
  {\it Uspehi Fiz. Nauk.} {\bf 9}, 141-166 (1929).  [in Russian]

\item\label{Khintchine MGU30}        %% 51 %%
 A.Ya.  Khintchine,
 Die Maxwell-Boltzmannsche Energieverteilung als Grenswertsatz
 der Wahrscheinlichkeitsrechung,
  {\it Trudy Sem. MGU Teor. Veroyatnost. i Mat. Statst.}
 {\bf 1}, 1-11 (1930). %!%

\item\label{Khintchine ATTUARI32}       %% 52 %%
 A.Ya.  Khintchine,
Sulle successioni stazionarie di eventi,
 {\it Giornale dell'Istituto Italiano degli Attuari},
{\bf 3}, No 3,  267-272 (1932).  %[?]
\\ {\footnotesize{[In  Kolmogorov  this references is quoted as [57]
 with {\bf 4}  (1933) 3]}}

\item\label{Khintchine MatSb32a}       %% 53 %%
 A.Ya.  Khintchine,
 Zur additiven  Zahlentheorie,
  {\it Mat. Sb.} {\bf 39}, No 3, 27-34 (1932). %!%

\item\label{Khintchine MatSb32a}    %% 54 %%
 A.Ya.  Khintchine,
Ueber eine Ungleichung,
  {\it Mat. Sb.} {\bf 39}, No 3, 35-39 (1932). %!%

\item\label{Khintchine MatSb32c}     %% 55 %%
 A.Ya.  Khintchine,
 Sur les classes d'\'ev\'enements \'equivalents,
  {\it Mat. Sb.} {\bf 39}, No 3, 40-43 (1932). %!%

\item\label{Khintchine MatSb32d}        %% 56 %%
 A.Ya.  Khintchine,
 Remarques sur les suites  d'\'ev\'enements ob\'eissants \'a la loi des
 grands nombres,
  {\it Mat. Sb.} {\bf 39}, 115-119 (1932).

\item\label{Khintchine MatSb32e}             %% 57 %%
 A.Ya.  Khintchine,
 The mathematical theory of stationary queues,
  {\it Mat. Sb.} {\bf 39}, No 4, 73-84 (1932) [in Russian].

\item\label{Khintchine MathAnn32b}             %% 58 %%
 A.Ya.  Khintchine,
Zur Birkoffs L\"osung der Ergodenproblems,
  {\it Math. Annalen} {\bf 107}, 485-488 (1932).

\item\label{Khintchine MatSb33a}                  %% 59 %%
 A.Ya.  Khintchine,
On the mean time of non-attendance of a machine,
  {\it Mat. Sb.} {\bf 40}, 119-123 (1933) [in Russian]. %!%

\item\label{Khintchine MatSb33b}              %% 60 %%
 A.Ya.  Khintchine,
Ueber stazion\"are Reihen zuf\"alliger Variablen
  {\it Mat. Sb.} {\bf 40}, 124-128 (1933).

\item\label{Khintchine MatSb33c}              %% 61 %%
 A.Ya.  Khintchine,
Ueber ein metrisches Problem der additiven Zahlentheorie
  {\it Mat. Sb.} {\bf 40}, 180-189 (1933).

\item\label{Khintchine ZAMM33}             %% 62 %%
 A.Ya.  Khintchine,
Zur mathematischen Begr\"undung des statistischen Mechanik,
  {\it Z. angew. Math. Mech.} {\bf 13}, 101-108 (1933).

\item\label{GelfandKhintchine MGU33}        %% 63 %%
A.O. Gel'fond and  A.Ya.  Khintchine, Gram's determinants for
stationary series,
  {\it Uchenue Zapiski Moscow Gov. Univ.} {\bf 1}, 3-5 (1933). %!%

\item\label{Khintchine USA33}                %% 64 %%
 A.Ya.  Khintchine,
The method of spectral reduction in classical dynamics,
  {\it Proc. Nat. Acad. Sci. U.S.A.} {\bf 19}, 567-573 (1933).

\item\label{Khintchine ASYMPROB33-36-48}    %% 65 %%
 A.Ya.  Khintchine, %%% (1948)
 {\it Asymptotische Gesetze der Wahrscheinlichkeits\-rechnung},
  Julius Springer, Berlin  (1933).%; ONTI, Moscow, 1936 [in Russian]
  %%  pp. 77.
Reprinted by Chelsea Publ., New-York,  1948. %% [X]
\\ {\footnotesize{
 There is another Edition (in Russian)
=(\it Asimptoticheskie zakony teorii veroyatnosfei)
{\it Asymptotic Laws  of the Theory of Probability}, O.N.T.I., Moscow (1936).}}

\item\label{Khintchine MGU34}              %% 66 %%
 A.Ya.  Khintchine,
Zur mathematischen Begr\"undung der Maxwell-Boltzmannschen
Energieverteilung,
  {\it Uchenue Zapiski Moscow Gov. Univ.} {\bf 2}, 35-38 (1934). %!%

\item\label{Khintchine MatAnn34-UMN38}     %% 67 %%
 A.Ya.  Khintchine,
Korrelationstheorie der station\"aren stochasti\-schen Prozesse
  {\it Math. Annalen} {\bf 109}, 604-615 (1934);
{\it Uspehi Mat. Nauk} {\bf 5}, 42-51 (1938).

\item\label{Khintchine ComMath34}         %% 68 %%
 A.Ya.  Khintchine,
Eine Versch\"arfung des Poincar\'eschen 'Wieder\-kehr\-satzes',
  {\it Compositio Math.} {\bf 1}, 177-179 (1934).

\item\label{Khintchine MatSb34a}           %% 69 %%
 A.Ya.  Khintchine,
Fourierkoeffizienten l\"angs Bahnen im Phasenraum,
  {\it Mat. Sb.} {\bf 41}, 14-16 (1934).

\item\label{Khintchine CHANCE34}             %% 70 %%
 A.Ya.  Khintchine,
{\it Chance and the Way Science Treats It},
 ONTI, Moscow, 1934 [in Russian].

\item\label{Khintchine MatSb34b}           %% 71 %%
 A.Ya.  Khintchine,
Eine arithmetische Eigenschaft der summierbaren Funktionen,
  {\it Mat. Sb.} {\bf 41}, 11-13 (1934).

\item\label{Khintchine PROBRussia34}           %% 72 %%
 A.Ya.  Khintchine,
The theory  of probability in the pre-revolutionary Russia and in
the Soviet Union,
  {\it Front Nauki i Techniki} {\bf 7}, 36-46 (1934) [in Russian].
  %!%

\item\label{Khintchine CONTFRACT35-49-63}    %% 73 %%
 A.Ya.  Khintchine, %%%
 {\it Continued Fractions},
ONTI, Moscow 1935; 2nd Edn GTTI, Moscow 1949; 3rd Edn Fizmatgiz,
Moscow, 1961; 4th Edn Nauka, 1978, 112 p. [in Russian]. English
Translation (by P. Wynn from the 3rd 1961 Russian edition),
Norddhoff, Groningen, 1963, pp.
101. %!%
%%% {[Xerocopied at MaPhySto, Jan 2002]}}

\item\label{Khintchine ATTUARI35}     %% 74 %%
 A.Ya.  Khintchine,
Sul dominio di attrazione della legge di Gauss,
 {\it Giornale dell'Istituto Italiano degli Attuari},
{\bf 6} No 4,  378-393 (1935).
%%%%%%%%%%
\\ {\footnotesize{Both Kolmogorov in [77] and Gnedenko in [79] cite this
paper in an wrong way, that is
Vol 7 (1936) pp. 3-18.
%%%%%%%
From the Obituary of Khintichine by Gnedenko(1960/61): For the
case of independent, identically distributed components Khintchine
(Khinchin)  succeeded, simultaneously with P. L\'evy and W. Feller
but independently from them, in finding the necessary and
sufficient conditions for convergence to the normal law [79 = this
paper] Papers [47 = in German on Math. Annalen, 1929 a] and [48 =
in German on Math. Annalen, 1929 b] deserve special mention
because they may be considered as initiating the current studies
of "large deviations".}}

\item\label{Khintchine ComMath35}              %% 75 %%
 A.Ya.  Khintchine,
Metrische Kettenbr\"uchprobleme,
  {\it Compositio Math.} {\bf 1}, 361-382 (1935).

\item\label{Khintchine MatAnn35}                 %% 76 %%
 A.Ya.  Khintchine,
Neuer Beweis und Verallgemeinerung eines Hurvitzschen Satzes,
  {\it Math. Annalen} {\bf 111}, 631-637 (1935).

\item\label{Khintchine ComMath36}          %% 77 %%
 A.Ya.  Khintchine,
Zur metrischen Kettenbr\"uchtheorie,
  {\it Compositio Math.} {\bf 3}, 276-285 (1936).

\item\label{Khintchine MatAnn36}             %% 78 %%
 A.Ya.  Khintchine,
Ein Satz \"uber lineare diophantische Approximationen,
  {\it Math. Annalen} {\bf 113}, 398-415 (1936).

\item\label{Khintchine ATTUARI36}             %% 79 %%
 A.Ya.  Khintchine,
Su una legge dei grandi numeri generalizzata,
 {\it Giornale dell'Istituto Italiano degli Attuari},
{\bf 7},  365-377 (1936).
%%%%%%%%%
\\{\footnotesize{[This is the ref [79] by Kolmogorov.
Gnedenko refers to it in [74], vol. 6 (1935), pp 371-393. The
different references to the Italian papers [74], [79] in this
journal  is a source of confusion that must be clarified: it is a
comedy of errors!]}}

\item\label{Khintchine UMN36}        %% 80 %%
 A.Ya.  Khintchine,
Metric problems in the theory of irrational numbers, {\it Uspehi
Mat. Nauk} {\bf 1}, 7-32 (1936).

\item\label{Khintchine BMGU37a}         %% 81 %%
 A.Ya.  Khintchine,
A new derivation of a formula of P. L\'evy {\it Bull. Moscow Gov.
Univ.} {\bf 1}, vyp. 1, 1-5 (1937) [in Russian]. %!%

\item\label{Khintchine BMGU37b}            %% 82 %%
 A.Ya.  Khintchine,
On the arithmetic of distribution laws, {\it Bull. Moscow Gov.
Univ.} {\bf 1},  vyp. 1, 6-17 (1937) [in Russian]. %!%

\item\label{Khintchine BMGU37c}             %% 83 %%
 A.Ya.  Khintchine,
On a property of characteristic functions, {\it Bull. Moscow Gov.
Univ.} {\bf 1},  vyp. 5, 1-3 (1937) [in Russian]. %!%

\item\label{Khintchine BMGU37d}           %% 84 %%
 A.Ya.  Khintchine,
Invariant classes of distribution laws {\it Bull. Moscow Gov.
Univ.} {\bf 1},  vyp. 5, 4-5 (1937) [in Russian]. %!%

\item\label{Khintchine BMGU37e}       %% 85 %%
 A.Ya.  Khintchine,
Examples of random variables satisfying stable probability laws,
{\it Bull. Moscow Gov. Univ.} {\bf 1}, vyp. 5, 6-9 (1937) [in
Russian]. %!%

\item\label{KhintchineLevy CRASP37}   %% 86 %%
 A. Ya. Khintchine and P.  L\'evy, Sur le lois stables,
  {\it C.R. Acad. Sci.  Paris} {\bf 202}  374-376 (1937). %%  [X]

\item\label{Khintchine BMGU37f}   %% 87 %%%
 A.Ya.  Khintchine,
Ueber die angen\"aherte Aufl\"osung linearer Gleichungen in ganzen
Zahlen, {\it Acta Arith.} {\bf 2}, 161-172 (1937).
%%%%%%%%%%%
\\ {\footnotesize{see  {\it Mat. Sb.} {\bf 32}, 203-219 (1925).}}

\item\label{Khintchine ComMath37}    %% 88 %%
 A.Ya.  Khintchine,
Ueber singul\"are Zahlensysteme,
  {\it Compositio Math.} {\bf 4}, 424-431 (1937).

\item\label{Khintchine BIMT37a}       %% 89 %%
 A.Ya.  Khintchine,
Absch\"atzungen beim Koneckerschen Approximationensatz, {\it Bull.
Inst. Math. Tomsk} {\bf 1}, 263-265 (1937).

\item\label{Khintchine BIMT37b}       %% 90 %%
 A.Ya.  Khintchine,
Ueber Klassnkonvergens von Verteilunsgesetzen, {\it Bull. Inst.
Math. Tomsk} {\bf 1}, 258-262 (1937).

\item\label{Khintchine MatSb37}     %% 91 %%%
 A.Ya.  Khintchine,
Zur Theorie der unbeschr\"anktteilbaren Verteilungs\-gesetze, {\it
Mat. Sbornik} [{\it Rec. Mat. [Mat. Sbornik}], New Series {\bf 2}
(44), No 1,  79-119 (1937). %%   [X]
%%%%%%%%

\item\label{Khintchine LIMITDISTR38}   %% 92 %%
 A.Ya. Khintchine, %% (1938)
{\it Limit Distributions for the Sum of Independent Random Variables}.
O.N.T.I., Moscow, 1938, pp. 115. [in Russian] %%  [X]
%%%%%%%%%%
%% \\ {\footnotesize{%%%%%%%%
%% From the Obituary of Khintichine by Gnedenko(1960/61): The theory
%% of summation inspired  Khintchine (Khinchin)   to write three
%% monographs. The first [35 = Fundamental Laws of the Theory of
%% Probability, in Russian, 2nd Edn 1932] was published in 1927 after
%% the completion of a special course on summation theory given by
%% Khintchine (Khinchin) at Moscow University. The second monograph
%% [65 = Asymptotic Laws  of the Theory of Probability, in German,
%% 1933 (in Russian? 1936)] connected   the classical problems of
%% summation with the theory of Markov proceses and with the then
%% recent investigations of Kolmogorov and Petrovsky.
%% The third monograph  [92 = this book!] gave a constructive development of
%% general limit theorems for sums of independent summands and their
% application  to the classical problem of convergence of normed
%% sums to the normal law. This book was also preceded by a special
%% course of lectures in Moscow University. This course attracted the
%% interest of A.A. Bobrov, D.A. Raikov, and myself to the summation theory.}}

\item\label{Khintchine MatSb38a}  %%% 93 %%%
 A.Ya. Khintchine, %% 1938 %%
Two theorems about stochastic processes with increments of the
same type, {\it Mat. Sb.} {\bf 3} [45], 577-584 (1938) [in
Russian] %%  [X].

\item\label{Khintchine MatSb38b}  %%% 94 %%%
 A.Ya. Khintchine,
 Zur Methode der willk\"{u}rlichen Funktionen,
{\it Mat. Sb.} {\bf 3} [45], 585-589 (1938). %% [X].

\item\label{Khintchine ITMI38}  %%% 95 %%%
 A.Ya. Khintchine,
 On unimodal distributions,
{\it Izv. Tomsk. Math. Inst.} {\bf 2}, 1-7 (1938) [in Russian].

\item\label{Khintchine IANS38}  %%% 96 %%%
 A.Ya. Khintchine,
 The theory of damped spontaneous effects,
{\it Izv. Akad. Nauk SSSR} (Ser. Mat.) {\bf 3}, 313-332  (1938)
[in Russian].

\item\label{Khintchine IRRNUMB38}  %%% 97 %%%
 A.Ya. Khintchine,
 Introduction to irrational numbers; material for the use of teachers,
{\it Narcompross RSFSR}, 9-12 (1938) [in Russian]. Reprinted in
{\it Mat. v. \v{S}kole} No 3, 32-34 (1939).

\item\label{KhintchineDorf SKOLE38}  %%% 98 %%%
 P.Ya. Dorf and  A.Ya. Khintchine,
 Complex numbers; material for the use of teachers,
 {\it Narcompross RSFSR}, 39-47 (1938) [in Russian]. %!%

\item\label{Khintchine MatSb39}  %%% 99 %%%
 A.Ya. Khintchine,
On the addition of sequences of natural numbers, {\it Mat. Sb.}
{\bf 4} (48), 161-166 (1939) [in Russian].

\item\label{Khintchine IANS30}  %%% 100 %%%
 A.Ya. Khintchine,
On local growth  of homogeneous stochastic processes without
future. {\it Izv. Akad. Nauk SSSR} (Ser. Mat.) {\bf 4}, 487-508
(1939) [in Russian].

\item\label{Khintchinef SKOLE39a}  %%% 101 %%%
 A.Ya. Khintchine,
 Fundamental concepts of mathematics in secondary schools,
{\it Mat. v. \v{S}kole} No 4, 4-22 (1939); No 5, 3-10 (1939).

\item\label{Khintchine SKOLE39b}  %%% 102 %%%
  A.Ya. Khintchine,
Many-sided realistic education of Soviet youth,
 {\it Mat. v. \v{S}kole} No 6, 1-7 (1939).

\item\label{Khintchine Molodaya39}  %%% 103 %%%
 A.Ya. Khintchine,
 On the teaching of mathematics,
{\it Molodaya Gvardia}  {No 9}, 142-150  (1939); reprinted in {\it
Mat. v. \v{S}kole} No 6, 1-7 (1939) under the title ``Many-sided
realistic education of Soviet youth''
[in Russian]. %!%
%% Included in the booklet %%% 160n %%%

\item\label{Khintchine SCHOOLS40}  %%% 104%%%
 A.Ya. Khintchine,
{\it Fundamental Mathematical Concepts and Notions  in Secondary
Schools} Uchpedciz, 1940, 51 pp. [in Russian].

\item\label{Khintchine UMN40}  %%% 105 %%%
 A.Ya. Khintchine,
 On the addition of sequences of natural numbers,
{\it Uspehi Mat. Nauk}  {\bf 7}?, 57-61  (1940) [in Russian].

\item\label{Khintchine SKOLE41a}  %%% 106 %%%
  A.Ya. Khintchine,
On  mathematical definitions in secondary schools,
 {\it Mat. v. \v{S}kole} No 1, 1-10 (1941). [in Russian]

\item\label{Khintchine SKOLE41b}  %%% 107 %%%
  A.Ya. Khintchine,
On  the concept of the ratio of two numbers,
 {\it Mat. v. \v{S}kole} No 2, 13-15 (1941).   [in Russian]

\item\label{Khintchine DOKL41a}  %%% 108 %%%
  A.Ya. Khintchine,
On  analytic methods of statistical mechanics, {\it Dokl. Akad.
Nauk SSSR} {\bf 33}, 438-441 (1941) [in Russian].

\item\label{Khintchine DOKL41b}  %%% 109 %%%
  A.Ya. Khintchine,
The mean value of summable functions in statistical mechanics,
 {\it Dokl. Akad. Nauk SSSR} {\bf 33}, 442-445 (1941) [in Russian].

\item\label{Khintchine DOKL41c}  %%% 110 %%%
  A.Ya. Khintchine,
On  intermolecular correlation,
 {\it Dokl. Akad. Nauk SSSR} {\bf 33}, 487-490 (1941) [in Russian].

\item\label{Khintchine DOKL42}  %%% 111 %%%
  A.Ya. Khintchine,
Laws of distribution of summable functions in statistical
mechanics,
 {\it Dokl. Akad. Nauk SSSR} {\bf 34}, 61-63 (1942) [in Russian].

\item\label{Khintchine MatSb43}  %%% 112 %%%
  A.Ya. Khintchine,
Sur un cas de corr\'elation a posteriori,
 {\it Mat. Sb.} {\bf 8} (54), 185-196 (1943) [in French].

\item\label{Khintchine STATMECH43-49}    %% 113 %%
 A.Ya. Khintchine, %% (1949)
{\it Mathematical Foundations of Statistical Mechanics}, GTTI,
Moscow, 1943, pp. 126  [in Russian].
 English Translation by G. Gamow  published by
 Dover, New-York, 1949, pp. 179. %% [X]

\item\label{Khintchine IANS43a}  %%% 114 %%%
 A.Ya. Khintchine,
On ergodic problem of statistical mechanics, {\it Izv. Akad. Nauk
SSSR}  {\bf 7}, 167-184  (1943) [in Russian].

\item\label{Khintchine IANS43b}  %%% 115 %%%
 A.Ya. Khintchine,
Convex functions an evolution theorems of statistical mechanics,
{\it Izv. Akad. Nauk SSSR}  {\bf 7}, 111-122  (1943) [in Russian].

\item\label{Khintchine MATHANAL43}  %%% 116 %%%
 A.Ya. Khintchine,
{\it Eight Lectures on Mathematical Analysis},
 GTTI, Moscow 1943, 3rd Edition 1948 [in Russian].
% There are also some later editions.

\item\label{Khintchine IANS46}  %%% 117 %%%
 A.Ya. Khintchine,
On a problem of Chebyshev. {\it Izv. Akad. Nauk SSSR} (Ser. Mat.)
{\bf 10}, 281-294 (1946) [in Russian].

\item\label{GnedenkoKhintchine PROB46}  %%% 118 %%%
 B.V. Gnedenko and A.Ya. Khintchine,
{\it Elementary Introduction to the Theory of Probability},
 Moscow, 1946, 4th Edition 1948, 5th Edition 1961, 6th Edition 1964 [in Russian].
%%%%%%%
\\ {\footnotesize{
Translations are published in Warsaw, 1952, 1954;
 Bucharest, 1953; Prague 1954; Budapest, 1954; Berlin, 1955; Paris,
 1960; Buenos Aires, 1960; New York, 1961.}}

\item\label{Khintchine PDAGOG46}  %%% 119 %%%
 A.Ya. Khintchine,
On formalism in the teaching of mathematics,
{\it Soviet Pedagogics}, No 11-12, 21--27 (1944); reprinted by
{\it Izv. Akad. Pedagog. Nauk RSPSR} (Ser. Mat.) {\bf 4}, 7-20 (1946). %!%
%% Included in the booklet %%% 160n %%%

\item\label{Khintchine NUMTHEORY47}  %%% 120 %%%
 A.Ya. Khintchine,
{\it Three Pearls of Number Theory}, GTTI, Moscow, 1947; 2nd Edn
1948; Nauka, Moscow, 3rd Edn 1979, pp. 64 [in Russian]; Ukrainian
Edn 1949; Berlin, 1950; Rochster, N.Y., 1952; Tokyo, 1956. %!%

\item\label{Khintchine IANS47}  %%% 121 %%%
 A.Ya. Khintchine,
Two theorems connected with  a problem of Chebyshev. {\it Izv.
Akad. Nauk SSSR} (Ser. Mat.)  {\bf 11}, 105-110 (1946) [in
Russian].

\item\label{Khintchine DOKL47a}  %%% 122 %%%
  A.Ya. Khintchine,
A limiting case of Kronecker's approximation theorem,
 {\it Dokl. Akad. Nauk SSSR} {\bf 56}, 563-565 (1947) [in Russian].

\item\label{Khintchine DOKL47b}  %%% 123 %%%
  A.Ya. Khintchine,
On a general theorem in the theory of Diophantine approximation,
 {\it Dokl. Akad. Nauk SSSR} {\bf 56}, 679-681 (1947) [in
 Russian]. %!%

\item\label{Khintchine DOKL48a}  %%% 124 %%%
  A.Ya. Khintchine,
A transient theorem for singular systems of linear equations,
 {\it Dokl. Akad. Nauk SSSR} {\bf 59}, 217-218 (1948) [in Russian].

\item\label{Khintchine DOKL48b}  %%% 125 %%%
  A.Ya. Khintchine,
On  the theory of linear Diophantine approximation,,
 {\it Dokl. Akad. Nauk SSSR} {\bf 59}, 865-867 (1948) [in Russian].

\item\label{Khintchine UMN48a}  %%% 126 %%%
  A.Ya. Khintchine,
Dirichlet's principle in Diophantine approximation,,
 {\it Uspehi Mat. Nauk} {\bf 3}, vyp. 3, 1-28 (1948) [in Russian].

\item\label{Khintchine IANS48a}  %%% 127 %%%
 A.Ya. Khintchine,
Quantitative concept of the approximation theory of Kronecker,
{\it Izv. Akad. Nauk SSSR} (Ser. Mat.)  {\bf 12}, 113-122 (1948) %!%
[in Russian].

\item\label{Khintchine UMN48b}  %%% 128 %%%
  A.Ya. Khintchine,
On some applications of the method of an additional variable, {\it
Uspehi Mat. Nauk} {\bf 3}, vyp. 6, 188-200 (1948) [in Russian]. %!%

\item\label{Khintchine IANS48b}  %%% 129 %%%
 A.Ya. Khintchine,
Regular systems of linear equations and the general problem of
Chebyshev, {\it Izv. Akad. Nauk SSSR} (Ser. Mat.)  {\bf 12},
249-258 (1948) [in Russian].

\item\label{Khintchine IANS49}  %%% 130 %%%
 A.Ya. Khintchine,
On fractional parts of a linear form, {\it Izv. Akad. Nauk SSSR}
(Ser. Mat.)  {\bf 13}, 3-8 (1949) [in Russian].

\item\label{Khintchine UMN49}  %%% 131 %%%
  A.Ya. Khintchine,
The simplest linear continuum,
 {\it Uspehi Mat. Nauk} {\bf 4}, vyp. 2,  180-197 (1949) [in Russian].

\item\label{Khintchine TMIANS}  %%% 132 %%%
  A.Ya. Khintchine,
On the analytical apparatus of physical statistics, %!%
 {\it Trudy Mat. Inst. Akad. Nauk SSSR} {\bf 33}, 1-56 (1950) [in Russian].

\item\label{Khintchine UMN50}  %%% 133 %%%
  A.Ya. Khintchine,
Statistical mechanics as a problem of the probability theory, %!%
 {\it Uspehi Mat. Nauk} {\bf 5}, vyp. 3, 3-46 (1950) [in Russian].

\item\label{Khintchine DOKL50}  %%% 134 %%%
  A.Ya. Khintchine,
On  sums of positive random variables,
 {\it Dokl. Akad. Nauk SSSR} {\bf 71}, 1037-1039 (1950) [in Russian].

\item\label{Khintchine UMZ50}  %%% 135 %%%
  A.Ya. Khintchine,
Limit theorems for  sums of positive random variables,
 {\it Ukrain. Mat. \u{Z}} {\bf 2}, 3-17 (1950) [in Russian].

\item\label{Khintchine QuantumStat51-56}  %%% 136 %%%
  A.Ya. Khintchine,
{\it The Mathemtical Theory of  Quantum Statistics} GTTI, Moscow
1951 [in Russian]; Berlin  1956 [in German].

\item\label{Khintchine DOKLAD51}  %%% 137 %%%
  A.Ya. Khintchine,
On the distribution of laws of 'occupation numbers' in  quantum
statistics, {\it Dokl. Akad. Nauk SSSR} {\bf 78}, 461-463 (1951)
[in Russian].

\item\label{Khintchine TMIANS50}  %%% 138 %%%
  A.Ya. Khintchine,
On some general theorems in statistical physics,
 {\it Trudy  Mat. Inst. Akad. Nauk SSSR} {\bf 38}, 345-365 (1951) [in Russian].

\item\label{Khintchine NUMTH51}  %%% 139 %%%
  A.Ya. Khintchine,
Elements of number theory, in {\it Encyclopedia of Elementary
Mathematics}
 GTTI, Moscow, {\bf 1}, 255-353 (1952). [in Russian] %!%

\item\label{Khintchine DOKLAD52}  %%% 140 %%%
  A.Ya. Khintchine,
On classes of equivalent events,
 {\it Dokl. Akad. Nauk SSSR} {\bf 85}, 713-714 (1952) [in Russian].

\item\label{Khintchine CONTPHYS52}  %%% 141 %%%
  A.Ya. Khintchine,
Method of arbitrary functions and the struggle against idealism in
the  probability theory, in the collection {\it Philosophical
Problems of Contemporary Physics} edited by {\it Akad. Nauk SSR},
1952, pp. 522-538 [in Russian]; French Translation: "Questions
Scientifiques V" , in {\it Editions de la Nouvelle Critique},
Paris, 1954, Vol. 1, pp. 7-24. %!%

\item\label{Khintchine SOVIETPROB52}  %%% 142 %%%
  A.Ya. Khintchine,
Soviet school of the theory of probability,
 {\it Chinese Math. J.} {\bf 1}, 1-7 (1952) [in Russian].

\item\label{Khintchine UMN53}  %%% 143 %%%
  A.Ya. Khintchine,
Notion of entropy in the  probability theory,
 {\it Uspehi Mat.  Nauk} {\bf 8}, vyp. 3, 3-20 (1953) [in Russian];
Bucharest, 1955; Berlin, 1956; U.S.A. 1957. %!%

\item\label{Khintchine ANALYSIS53-57}  %%% 144 %%%
  A.Ya. Khintchine,
{\it A Short Course of Mathematical Analysis},
 GTTI, Moscow 1953; 3rd Edn 1957 [in Russian]; Bucharest, 1955; Beijing, 1957.

\item\label{AlexandrovKhintchine Kolmoorov53}  %%% 145 %%%
  P.S. Alexandrov and A.Ya. Khintchine,
Andrei Nikolaevich Kolmogorov (on his fiftieth anniversary),
 {\it Uspehi Mat.  Nauk} {\bf 8}, vyp. 3, 177-200 (1953) [in
 Russian]; %!%

\item\label{Khintchine TRMIANS55}  %%% 146 %%%
  A.Ya. Khintchine,
Mathematical methods in the theory of mass services (queueing
theory),
 {\it Trudy  Mat. Inst. Akad. Nauk SSSR} {\bf 49}, 1-123 (1955) [in Russian];
 Beijing, 1958. %!%

\item\label{Khintchine IANS55}  %%% 147 %%%
 A.Ya. Khintchine,
Symmetric functions over multidimensional surfaces, A contribution
to the memorial volume to Alexander Alexandrovich Andronov, {\it
Izdat. Akad. Nauk SSSR}, 541-574  (1955) [in Russian].

\item\label{Khintchine TVP56a}         %% 148 %%
 A.Ya. Khintchine, %% (1956)
 Streams of random events without future,
{\it Teoriya Veroyatnostei i  Primenenie} ({\it Probability Theory
and Applications}) {\bf 1} No 1, 3-18 (1956). %!%

\item\label{Khintchine TVP56b}         %% 149 %%
 A.Ya. Khintchine, %% (1956)
 On Poisson streams of random events,
{\it Teoriya Veroyatnostei i  Primenenie} ({\it Probability Theory
and Applications})
 {\bf 1} No 3, 320-327 (1956).

\item\label{Khintchine UMN56}         %% 150 %%
 A.Ya. Khintchine,
 On fundamental theorems of information theory,
 {\it Uspehi Mat. Nauk} {\bf 11}, vyp. 1, 17-75 (1956). %!%
\\
{\footnotesize{The English translation of this paper has appeared
 in {\it Mathematical Foundations of Information Theory},
 Dover, New-York, 1957.}}

\item\label{KhintchineYaglom CHILDREN59}         %% 151 %%
 A.Ya. Khintchine and A.M. Yaglom,
 The first acquaintance with probability theory,
 {\it Children's Encyclopedia}, edited by {\it Akad. Pedagog. Nauk RSFSR},
{\bf 3}, 211-220 (1959).
\\
{\footnotesize{With this joint article devoted to the education of
children towards probability the Bibliography of A. Ya Khintchine
by Gnedenko ends.}}

\vspace{0.25truecm} \hrule \vspace{0.25truecm}

%\item\label{Khintchine QUEUEING55-60} %%% 152n=146 %%%
% A.Ya. Khintchine,
% {\it Mathematical Methods in the Theory of Queueing},
%   Charles Griffin, London 1960, pp. 120.
% [Translated from the 1955 Russian edition]   [X]
% \\ {\footnotesize{[No 7  of the Griffin's Statistical Monographs
% \& Courses edited by M.G. Kendall.]
% [Russian edition {\it Mathematicheskie metody teorii massovogo
% obsluzhivania} in Trudy Matematicheskogo Instituta
%im. V.A. Steklova, Vol. 49, Izd AN SSR, 1955]}}
%%% Avalable  MATH-Library FU-BERLIN but XEROCOPIED with REVISION
%% from Gorenflo's personal book %%%%

%%%
\item\label{Khintchine PAPERSQUEUEING63} %%% 153n %%%
 A.Ya. Khintchine,
 {\it Papers on the Mathematical Theory  of Queues}
({\it Raboty po matematicheskoi teorii massovogo obsluzhivaniya}),
Fizmatgz, Moskva, 1963.

%\item\label{Khintchine MATHANAL60}
% A.Ya.  Khintchine (A. Khinchin), %%% %%% 154n=144 %%%
% {\it A Course of Mathematical Analysis}
% (Industan Publ. Co, Delhi,  1960), pp. 668.
%[Translated from the 3rd 1957 Russian edition (1st ed. 1953, 2nd
%ed. 1954).
%\\ {\footnotesize{[Partially xerocopied at MaPhySto, Jan 2002:
%Contents + Conclusions (Short historical sketch)]}}

   \item\label{Khintchine INFTHEORY57} %%% 155n %%%
 A.Ya. Khintchine, %% (1957)
 {\it Mathematical Foundations of Information Theory},
 Dover, New-York, 1957, pp. 120 [Translated from the Russian]  %% [X]
\\ {\footnotesize{This Dover edition is a translation of
the the following two papers by Khintchine:
\\
1) {\it The Entropy Concept in Probability Theory}, pp. 1-28: see
[143] in Gnedenko list, published in
 {\it Uspehi Mat. Nauk} {\bf 8}, 3-20 (1953).
\\
2) {\it On the Fundamental Theorems of Information Theory}, pp.
29-120;
 see [150] in Gnedenko list, published in
 {\it Uspehi Mat. Nauk} {\bf 11}, 17-75 (1956).
\\
In the first paper the Author develops the concept of entropy in
probability theory as a measure of uncertainty of a finite
'scheme', and discusses a simple application to coding theory. The
second paper investigates the restrictions previously placed on
the study of sources, channels and codes and attempts, "to give a
complete, detailed proof of both... Shannon theorems, assuming any
ergodic source and any stationary channel with a finite memory."
Besides surveying existing literature, the Author includes many
original contributions not readily available elsewhere.}}
%%%%%%%%%%%%%%%%%%%%%%%%%%%%%%%%%%%%%

%\item\label{Khintchine MATHANAL58} %%% 156n=144 %%%
% A.Ya.  Khintchine, %%%
% {\it A Short Course of Mathematical Analysis}
% (Vilnius, Lithuania,  1958), pp. 579 (in Lithuanian).

% \item\label{Khintchine FRACTIONS61}
% A.Ya.  Khintchine, %%%
% {\it Continued Fractions}
% (Phys.-Mat., Moscow,  3rd Edition, 1961), pp. 112 (in Russian).
% It is already cited.

 \item\label{Khintchine vonMises61a} %%% 157n %%%
 A.Ya.  Khintchine, %%%
 R.~von Mises' frequency theory and modern ideas of Probability Theory,
 {\it  Questions of Philosophy}, No. 1, 91--102 (1961)
  (in Russian).

  \item\label{Khintchine vonMises61b} %%% 158n %%%
 A.Ya.  Khintchine, %%%
 R.~von Mises' frequency theory and modern ideas of Probability Theory,
 {\it  Questions of Philosophy}, No. 2, 77--89 (1961)
  (in Russian).

\item\label{Khintchine TVP62}         %% 159n %%
 A.Ya. Khintchine, %% (1962)
 On Erlang's formulas in the queueing theory,
{\it Teoriya Veroyatnostei i  Primenenie} ({\it Probability Theory
and Applications})
 {\bf 7} No 3, 330-335 (1962).

\item\label{Khintchine Educational}
A.Ya. Khintchine, %% (1961-1962)
On educational effect of the mathematical classes, In {\it
Mathematical Education}, No 6, 7--28 (1961); reprinted by {\it
Matematika v \v{S}kole}, No 3, 30--44 (1962).
% included into collection of papers %%%160n%%%

\item\label{Khintchine Reasoning}
A.Ya. Khintchine, %% (1961)
on so called ``problems on reasoning'' in the course of
arithmetics, In {\it Mathematical Education}, No 6, 29--36 (1961).
% included into collection of papers %%%160n%%%
% Comment by B.V. Gnedenko in \cite[p. 161]{Khintchine Pedagogical}:
% As I know this article was never published inter vivos of A.Ya. Khintchine.
% I have found it prepared (written by pencil in a school notebook). Since
% these questions were interested him in 1938-1939, almost sure this article
% was was written at that time.

 \item\label{Khintchine Pedagogical} %% %%% 160n %%%
 A.Ya.  Khintchine, %%%
 {\it Pedagogical Papers}
 ({\it Pedagogicheskie stat'i}), Prosvestchenie, Moscow,   1963, pp. 204 (in Russian).

\end{list}

\end{document}